\documentclass[12pt,a4paper]{article}
\usepackage[T1]{fontenc}
\usepackage{amssymb}
\usepackage{amstext} 
\usepackage{amsmath} 
\usepackage{amscd}
\usepackage{hyperref}
\usepackage[all]{xypic}
\usepackage{rotating}
\usepackage{graphicx}
\usepackage{color}
\usepackage{url} 
\usepackage[french]{babel}
\usepackage{pdfcomment}
\usepackage{theorem}
\usepackage{mathrsfs}
\usepackage{ulem}

\newlength{\secskip}
\newlength{\thsecskip}
\newlength{\ssecskip}
\newlength{\thssecskip}
\newlength{\sssecskip}
\newlength{\thsssecskip}
\makeatletter
\renewcommand\section{\@startsection{section}{1}{\z@}%
                                   {\secskip}
                                   {2.3ex \@plus.2ex}
                                   {\normalfont\Large\bfseries}}
\renewcommand\subsection{
			  \@startsection{subsection}{2}{\z@}%
                                     {\ssecskip}
                                     {1.5ex \@plus .2ex}
									 {\normalfont\normalsize\bfseries}}
\renewcommand\subsubsection{
              \@startsection{subsubsection}{3}{\z@}%
                                     {\sssecskip}
                                     {\z@}
									 {\normalfont\normalsize\bfseries}}
\makeatother

\makeatletter
\renewcommand\l@subsection{\@dottedtocline{2}{1.5em}{2.7em}}
\makeatother

\makeatletter
\renewcommand\l@subsubsection{\@dottedtocline{3}{3.8em}{3.5em}}
\makeatother
\newcommand\rref[1]{{\rm\ref{#1}}}
\theoremstyle{change}
\theorembodyfont{\slshape}
\theorembodyfont{\rmfamily}

\newtheorem{fcond}[subsection]{Condition}
\newtheorem{fscond}[subsubsection]{Condition}
\newtheorem{fsscond}[paragraph]{Condition}

\newcommand{\Spec}{\mathop{\mathrm{Spec}}}

\newcommand{\dem}{{\sl \noindent D\'emonstration\hskip0.1em: }} 
\newcommand{{\rH}}{\mathrm{H}}
\newcommand{\Hom}{{\rm Hom}}

\newcommand{\id}{{\mathrm{Id}}}

\newcommand{\im}{\mathrm{Im}}

\newcommand{\Top}{\mathrm{top}}

\newcommand{\abs}{\mathrm{abs}}

\newcommand{\FF}{{\mathbb F}}
\newcommand{\CC}{{\mathbb C}}
\newcommand{\RR}{{\mathbb R}}
\newcommand{\ZZ}{{\mathbb Z}}
\newcommand{\QQ}{{\mathbb Q}}
\newcommand{\NN}{{\mathbb N}}

\newcommand{\PP}{{\mathbb P}}

\newcommand{\GG}{{\mathbb G}}

\newcommand{\soul}[1]{\underline{#1}}
\newcommand{\Aa}{\mathbb{A}} 
\newcommand{\ade}{\mathbf{A}} 
\newcommand{\ide}{\mathbf{I}} 

\newcommand{\carrenoir}{\vrule height2.5mm width2mm depth 0mm}
\def\qed{ \hbox to5mm{}\hfill\llap{\carrenoir}}

\renewcommand{\setminus}{\smallsetminus}

\def\gG{\mathfrak{G}}

\def\gX{\mathfrak{X}}

\def\gY{\mathfrak{Y}}

\def\FI{\mathfrak{I}}

\def\Fm{\mathfrak{m}}

\newcommand\FX{\mathfrak{X}}
\newcommand\FZ{\mathfrak{Z}}
\newcommand\FY{\mathfrak{Y}}
\newcommand\FU{\mathfrak{U}}

\newcommand{\cB}{\mathscr{B}}
\newcommand{\cO}{\mathscr{O}}

\newcommand{\ol}{\overline}

\newcommand{\fl}{\rightarrow} 
\newcommand{\ffl}{\longrightarrow} 
\newcommand{\inj}{\hookrightarrow}

\newcommand{\sfl}[1]{\mathop{\fl}\limits^{#1}}

\newcommand{\flis}{\sfl{\sim}}

\newcommand{\surj}{\rlap{$\rightarrow$}\kern-1pt\rightarrow}
\newcommand{\ssurj}{\rlap{$\longrightarrow$}\kern-1pt\longrightarrow}
\newcommand{\surjgauche}{\rlap{$\leftarrow$}\kern-1pt\leftarrow}
\newcommand{\ssurjgauche}%
{\rlap{$\longleftarrow$}\kern-1pt\longleftarrow}

\newcommand{\mapdown}[1]
{\big\downarrow\rlap{$\vcenter{\hbox{$\scriptstyle#1$}}$}}

\newcommand{\mapup}[1]
{\big\uparrow\rlap{$\vcenter{\hbox{$\scriptstyle#1$}}$}}

\newcommand{\Mapdown}[1]
{\Big\downarrow\rlap{$\vcenter{\hbox{$\scriptstyle#1$}}$}}

\newcommand{\varfl}[1]{%
\setbox0=\hbox{$\;\;{\scriptstyle#1}\;\;\;$}%
\setbox1=\hbox to\wd0{$\;$\rightarrowfill$\;$}%
{\mathop{\box1}\limits^{\box0}}%
}
\newcommand{\varflspecial}[2]{%
\setbox0=\hbox{$\;\;{\scriptstyle#1}\;\;\;$}%
\setbox1=\hbox to\wd0{$\;$\rightarrowfill$\;$}
{\mathop{\box1}\limits^{#2}}%
}
\newcommand{\varflman}[2]{
\setbox0= \hbox to #1{$\;$\rightarrowfill$\;$}
{\mathop{\box0}\limits^{#2}}
}
\newcommand{\varmapstoman}[2]{%
\setbox0= \hbox to #1{$\;\mapstochar$\rightarrowfill$\;$}
{\mathop{\box0}\limits^{#2}}
}

\newcommand{\pr}{\mathrm{pr}}

\newcounter{nc}
\renewcommand{\thenc}{{\rm(\roman{nc})}}
\newenvironment{romlist}%
{\begin{list}{\thenc}{
\usecounter{nc} 
\parsep=0pt
\setlength  \labelwidth{\leftmargin}
\addtolength\labelwidth{-\labelsep}
}
}{\end{list}}
\newcounter{nnc}
\renewcommand{\thennc}{{\rm(\alph{nnc})}}
{\begin{list}{\thennc}{
\usecounter{nnc}
\parsep=0pt
\setlength  \labelwidth{\leftmargin}
\addtolength\labelwidth{-\labelsep}
}
}{\end{list}} 
\newcommand{\pauseromlist}%
{\global\edef\savecount{\arabic{nc}}\end{romlist}}
\newcommand{\finpauseromlist}%
{\begin{romlist}\setcounter{nc}{\savecount}}

\newcounter{ctnum}
\renewcommand{\thectnum}{\textup{(\arabic{ctnum})}}
\newenvironment{numlist}%
{\begin{list}{\thectnum}{
\usecounter{ctnum} 
\parsep=0pt
\leftmargin=0pt%
\setlength{\itemindent}{\labelwidth}%
\addtolength{\itemindent}{\labelsep}%
}
}{\end{list}}

\newcommand{\red}{{\mathrm{red}}}
\newcommand{\wh}[1]{\widehat{#1}}

\def\lgr{\longrightarrow}

\def\simlgr{\buildrel\sim\over\lgr}

\title{Fibr\'es principaux et ad\`eles }
\author{Philippe  Gille\\
Laurent  Moret-Bailly\smallskip\\
\small{\`A para\^{\i}tre dans \textsl{Annali della Scuola Normale Superiore di Pisa}}\\
\small{(accept\'e le 6 juin 2022)}}

\begin{document}

\date{}
\maketitle

\medskip

\noindent{\bf Abstract:}  We investigate topological properties of torsors in 
algebraic geometry over adelic rings.

\medskip

\noindent{\bf R\'esum\'e.} Nous  \'etudions les propri\'et\'es topologiques
des torseurs en g\'eom\'etrie alg\'ebrique sur des anneaux d'ad\`eles.

\medskip


\medskip


\setcounter{tocdepth}{2}


\bigskip

\addtocounter{section}{0}

\section{Introduction} \label{Sec_Intro}

L'article \cite{GGMB} porte sur les propri\'et\'es topologiques
des torseurs sous un  groupe alg\'ebrique sur des corps valu\'es, notamment
de caract\'eristique positive.
\'Etant donn\'e un corps global $K$, son anneau des ad\`eles $\ade_{{K}}$
est un anneau topologique et il est naturel de consid\'erer des questions analogues dans 
ce cadre. De fa\c con plus pr\'ecise, si $f: X\to Y$ est un morphisme
de $\ade_{{K}}$-sch\'emas qui est un $G$-torseur o\`u 
$G$ d\'esigne un $\ade_{{K}}$-sch\'ema en groupes (tous  de pr\'esentation finie), 
que peut-on dire de l'application 
$f_{ \Top}: X(\ade_{{K}})_\Top \to
X(\ade_{{K}})_\Top$?
Ici  $X(\ade_{{K}})_\Top$ (resp. $Y(\ade_{{K}})_\Top)$
d\'esigne l'espace
topologique associ\'e \`a $X$ (resp.  \`a $Y$)
selon Weil \cite{W}. En particulier, a-t-on des conditions naturelles pour garantir 
que l'image de $f_{\Top}$ est ferm\'ee ou localement ferm\'ee?

Si $G$ est affine lisse \`a fibres connexes, nous montrons  
 que l'image $I$ de $f_{\Top}$ est ferm\'ee  et que  
 l'application $X(\ade_{{K}})_\Top \to I$ est en outre une  
 $G(\ade_{\soul{K}})_\Top$-fibration principale 
 (cor.\ 6.3).
\smallskip
De fa\c con analogue \`a l'article \cite{GGMB}, notre approche ne n\'ecessite pas
d'imposer d'hypoth\`ese de compacit\'e locale pour les compl\'etions de sorte 
que certains de nos r\'esultats valent dans un cadre bien plus \'etendu qui inclut notamment
celui du corps de fonctions d'une courbe alg\'ebrique d\'efinie sur un corps arbitraire,
voir notamment le th\'eor\`eme 6.8. En fait, pour beacoup de nos \'enonc\'es, le corps 
\og arithm\'etique\fg\ $K$ sous-jacent ne joue aucun r\^ole, et les anneaux d'ad\`eles consid\'er\'es sont des produits restreints de familles plus ou moins arbitraires de corps valu\'es.

\smallskip

Passons en revue les diff\'erentes sections de l'article.
Le \S  \ref{SecTop} donne lieu \`a quelques sorites sur les topologies 
produit et le \S \ref{secLoc} revient sur la topologie \emph{forte} 
sur $X(A)$ pour $X$ un 
sch\'ema localement de type fini au dessus d'un anneau 
topologique $A$.
Le cas d'un anneau de valuation hens\'elien est particuli\`erement important et nous raffinons
plusieurs \'enonc\'es de \cite[\S 3]{GGMB}.
Le \S \ref{SecProdRestr} met en place le cadre ad\'elique dont nous avons besoin
et les topologies induites sur les points ad\'eliques des sch\'emas.
La section \ref{secPP} est le c{\oe}ur de l'article o\`u est d\'efinie la propri\'et\'e 
de \emph{presque propret\'e} pour un morphisme $f: X \to Y$ de
sch\'emas ad\'eliques, inspir\'ee par un th\'eor\`eme d'Oesterl\'e \cite[I.3.6]{O},
et qui est li\'ee \`a la propri\'et\'e d'image ouverte. C'est cette propri\'et\'e
que l'on \'etudie au \S \ref{sec:torseurs} dans le cadre des torseurs et
qui nous conduit aux r\'esultats principaux.
\medskip

{\sl Les auteurs remercient le rapporteur pour ses remarques.}

\section{Rappels topologiques}\label{SecTop}

\noindent{\bf D\'efinition 2.1.} 
{\it Soient $X$ et $Y$ deux espaces topologiques et $f:X\to Y$ une application continue. 
\begin{numlist}
\item Nous dirons que $f$ est un \emph{plongement topologique} si $f$ induit un hom\'eomorphisme de $X$ sur $f(X)$.
\item Nous dirons que $f$ est une \emph{fibration (topologique) triviale} s'il existe un espace $F$ et un diagramme commutatif
$$\xymatrix{X\ar[rr]^{\sim}_{u} \ar[dr]_{f} && Y\times F \ar[dl]^{\mathrm{pr}_{1}}\\
& Y&
}$$
o\`u $u$ est un hom\'eomorphisme.
\item Nous dirons que $f$ est une \emph{fibration topologique} s'il existe un recouvrement de $Y$ par des ouverts $U$ tels que les applications induites $f^{-1}(U)\to U$ soient des fibrations triviales.
\item Nous dirons que $f$ admet \emph{assez de sections locales} si, pour tout $x\in X$, il existe un voisinage $V$ de $f(x)$ dans $Y$ 
et une application continue $s:V\to X$ v\'erifiant $s(f(x))=x$ et  $f\circ s=\id_{V}$. Si, de plus, on peut prendre $V=Y$ pour tout $x$, nous dirons que $f$ admet \emph{assez de sections}.
\end{numlist}
}\smallskip
Une application ayant assez de sections locales est en particulier ouverte.

\subsection{Topologies produits }\label{SsecTopProduit} 

\noindent{\bf Lemme 2.2.}  
{\it Soient $V$ un ensemble et $(X_{v})_{v\in V}$ une famille index\'ee par
$V$ d'espaces topologiques; pour chaque $v\in V$, soit 
$Z_{v}$ un sous-espace de $X_{v}$. On pose $X=\prod_{v\in V}X_{v}$ et 
$Z=\prod_{v\in V}Z_{v}$, munis des topologies produits;
on note $\ol{Z_{v}}$ l'adh\'erence de 
$Z_v$ dans $X_v$, et $\ol{Z}$ celle de $Z$ dans $X$.

\begin{romlist}
\item\label{lemTopProduit5} Si $Z_{v}$ est ouvert dans $X_{v}$ pour tout $v$,
et si $Z_{v}=X_{v}$
pour presque tout $v$, alors $Z$ est ouvert dans $X$.
\item\label{lemTopProduit1} L'inclusion $Z\subset X$ est un plongement topologique, 
et $\ol{Z}=\prod_{v\in V}\ol{Z_{v}}$.
\item\label{lemTopProduit3} Si $Z_{v}$ est ferm\'e dans $X_{v}$ pour tout $v$, alors 
$Z$ est ferm\'e dans $X$.
\item\label{lemTopProduit4} Si $Z_{v}$ est localement ferm\'e dans $X_{v}$ pour tout
$v$, et ferm\'e pour presque tout $v$, 
alors $Z$ est localement ferm\'e dans $X$.
\qed
\end{romlist}
}

\smallskip

\dem \ref{lemTopProduit5} est clair par d\'efinition de la topologie produit.
Pour \ref{lemTopProduit1} (qui implique
trivialement \ref{lemTopProduit3}), voir \cite[I, \S 4, \no2, cor. de la prop. 3]{BTG}
et \cite[I, \S 4, \no3, prop. 7]{BTG}. 
Enfin, l'assertion
\ref{lemTopProduit4} r\'esulte de  \ref{lemTopProduit5} et \ref{lemTopProduit1}
appliqu\'ees en rempla\c{c}ant 
les $X_{v}$ par les $\ol{Z_v}$.\qed\medskip

\noindent La preuve du corollaire suivant est laiss\'ee au lecteur.\medskip

\noindent{\bf Corollaire 2.3.} 
{\it Soient $V$ un ensemble et $(f_{v}:X_{v}\to Y_{v})_{v\in V}$ une famille index\'ee par $V$ 
d'applications continues. Notons $f:X=\prod_{v\in V}X_{v}\to Y=\prod_{v\in V}Y_{v}$ l'application produit.
\begin{romlist}
\item\label{corTopProduit1} L'application naturelle $\im\,(f)\to\prod_{v\in V}\im\,(f_{v})$ est un hom\'eomorphisme.
\item\label{corTopProduit2} Si, pour tout $v\in V$, $f_{v}$ admet assez de sections, alors il en est de m\^eme de $f$.
\item\label{corTopProduit3} Si $f_{v}$ admet assez de sections locales pour tout $v\in V$, et admet assez de sections 
pour presque tout $v$, alors $f$ admet assez de sections locales.\qed
\end{romlist}
}

\section{Rappels locaux}\label{secLoc}
\subsection{Anneaux locaux topologiques}\label{ssecAnnLocTop}
Nous appellerons \emph{anneau local topologique} un anneau local 
$(A,\Fm)$, muni d'une topologie compatible avec sa structure d'anneau, 
et v\'erifiant de plus les propri\'et\'es suivantes:
\begin{itemize}
\item[(LT 1)]\label{CondLT1} $A^\times$ est ouvert dans $A$
(de fa\c{c}on  \'equivalente, $\Fm$ est ferm\'e);
\item[(LT 2)]\label{CondLT2} l'application 
$z\mapsto z^{-1}$ est continue sur $A^\times$ (de sorte que $(A^\times,.)$ 
est un groupe topologique).
\end{itemize}

\medskip

\noindent{\bf Exemples 3.1.} 
\begin{numlist}
\item\label{ExAnnLocTop1} Tout corps topologique (s\'epar\'e) est un
anneau local topologique.
\item\label{ExAnnLocTop2} Soit $(A,\Fm)$ un anneau local 
s\'epar\'e et \emph{lin\'eairement topologis\'e}, c'est-\`a-dire admettant 
une base $\FI$ de voisinages de $0$ form\'ee d'id\'eaux\footnote{Il faut prendre garde
que, en g\'en\'eral, les id\'eaux ne sont pas ouverts. 
Exemple: l'id\'eal nul, ou $(x)$ dans $k[[x,y]]$)}.  
$A$ est un anneau local topologique.
Voici les deux cas particuliers les plus importants:
\begin{romlist}
\item\label{ExAnnLocTop21} $A$ est un anneau de valuation (voir \S \ref{ssecVal}), muni de la topologie associ\'ee; on peut 
prendre pour $\FI$  l'ensemble des id\'eaux non nuls de $A$, ou encore l'ensemble des id\'eaux principaux non nuls;
\item\label{ExAnnLocTop22} $(A,\Fm)$ est un anneau local noeth\'erien, et  $\FI=\left\{\Fm^n\right\}_{n\geq1}$.
\end{romlist}
\end{numlist}

En pratique, nous utiliserons dans cet article les cas \ref{ExAnnLocTop1} et \ref{ExAnnLocTop2}\,\ref{ExAnnLocTop21}.

\subsection{Sch\'emas sur un anneau local topologique: topologie sur les sections. }\label{ssecTopLoc}
Soit $A$ un anneau local topologique. D'apr\`es \cite[Proposition 3.1]{ConAdelic}, pour tout $A$-sch\'ema $\FX$ 
localement de type fini on peut d\'efinir une topologie sur $\FX(A)$, donnant naissance \`a un foncteur $\FX\mapsto\FX(A)_{\Top}$
\`a valeurs dans les espaces topologiques, avec les propri\'et\'es suivantes:
\begin{itemize}
\item le foncteur $\FX\mapsto\FX(A)_{\Top}$ commute aux produits fibr\'es;
\item il transforme les immersions ferm\'ees (resp. ouvertes) en plongements topologiques ferm\'es (resp. ouverts);
\item la topologie de $\Aa^1(A)_{\Top}$ est la topologie donn\'ee sur $A$.
\end{itemize}
La topologie sur $\FX(A)$ est engendr\'ee par les ensembles de la forme
$$\left\{x\in \FU(A) \mid \varphi(x)\in V\right\}$$
o\`u $\FU$ d\'esigne un sous-sch\'ema ouvert de $\FX$, $V$ un ouvert de $A$ et $\varphi\in \Gamma(\FU,\cO_{\FX})$. En outre, on 
peut restreindre le choix de $\FU$ aux ouverts d'un recouvrement affine donn\'e, celui de $V$ \`a une base
d'ouverts de $A$, et celui de $\varphi$ \`a une fa\-mille de 
g\'en\'erateurs de la $A$-alg\`ebre $\Gamma(\FU,\cO_{\FX})$. Noter que si $\FX=\Spec(R)$ est affine, ceci peut s'exprimer 
en disant que l'inclusion naturelle de $\FX(A)=\Hom_{A\text{-alg}}(R,A)$ dans le produit $A^R$ est un plongement topologique.

\medskip

\noindent{\bf Remarque 3.2.}
Dans cette construction, le fait que $A$ soit local assure que si $(\FU_{i})_{i\in I}$ est un recouvrement ouvert 
de $\FX$, alors $\FX(A)$ est r\'eunion des  $\FU_{i}(A)$. Les conditions (LT~1) et (LT~2) de \ref{ssecAnnLocTop} 
entra\^{\i}nent que si $\FX=\Spec(R)$ est affine, et si $f\in R$, alors l'inclusion $D(f)\inj \FX$ est un plongement topologique 
ouvert pour les topologies \og affines\fg\ d\'ecrites ci-dessus.

\subsection{Le cas lin\'eairement topologis\'e. }\label{ssecLinTopLoc}
Soit $A$ un anneau local topologique. Supposons $A$ s\'epar\'e et lin\'eairement topologis\'e, et soit  $\FI$ une base
de voisinages de $0$ form\'ee d'id\'eaux stricts. 

Si $\FX$ est un $A$-sch\'ema localement de type fini, on peut alors construire $\FX(A)_{\Top}$  plus directement
\cite[\S 3.3]{GGMB} en d\'eclarant qu'une base de voisinages ouverts (et ferm\'es) d'un point $x\in\FX(A)$ est donn\'ee 
par les \og boules\fg
$$B_{\FX}(x,J):=\left\{y\in\FX(A)\mid y\equiv x \bmod J\right\}\qquad (J\in\FI)$$
o\`u \og$y\equiv x \bmod J$\fg\ signifie que $x$ et $y$ ont la m\^eme image dans $\FX(A/J)$. En d'autres termes, on
munit $\FX(A)$ de la topologie induite par l'application naturelle (injective, puisque $A$ est s\'epar\'e) 
$\FX(A)\inj \varprojlim_{J\in\FI}\FX(A/J)$, les ensembles $\FX(A/J)$ \'etant munis de la topologie discr\`ete.

Noter qu'en particulier $\FX(A)_{\Top}$ est toujours s\'epar\'e; en outre, pour tout id\'eal ouvert $J$ fix\'e, il
peut se d\'ecrire comme une somme topologique
$$
\FX(A)_{\Top}=\coprod_{x\in\Sigma}B_{\FX}(x,J)
\eqno{(3.1)}
$$
o\`u  $\Sigma\subset\FX(A)$ d\'esigne un syst\`eme repr\'esentatif des boules de rayon $J$.

\medskip

\noindent{\bf Proposition 3.3.} 
\textup{(Structure des morphismes lisses).} 
{\it Soit $(A,\Fm)$ un anneau local topologique s\'epar\'e et  \emph{lin\'e\-ai\-re\-ment topologis\'e}.
On suppose de plus que $A$ est \emph{hens\'elien}, et l'on fixe un id\'eal ouvert strict $J$ de $A$.

Soient $\FX$ et $\FY$ deux $A$-sch\'emas localement de type fini, et $f:\FX\to\FY$ un $A$-morphisme  \emph{lisse} de dimension relative $d$. 
On consid\`ere l'application continue induite $f_{\Top}:\FX(A)_{\Top}\to\FY(A)_{\Top}$.
\begin{numlist}
\item\label{lem-implicite1} \textup{(Fonctions implicites)} Soit $x \in \gX(A)$ et posons $y= f(x)$. Alors l'application 
$$B_{\FX}(x,J)\to B_{\FY}(y,J)$$ 
induite par $f_{\Top}$ est une fibration topologique triviale de fibre $J^d$%
\footnote{Ici, $J^d$ d\'esigne \'evidemment l'espace produit, et non l'id\'eal produit.}. 
En particulier, $f_{\Top}$ est ouverte, et est un hom\'eomorphisme local si $f$ est \'etale.
\item \label{lem-implicite2} L'image de $f_{\Top}$ s'\'ecrit de fa\c{c}on  unique 
$$ \im(f_{\Top})=\coprod_{\lambda\in L}B_{\lambda}$$
o\`u les $B_{\lambda}$ sont des  boules de rayon $J$ dans $\FY(A)_{\Top}$; en particulier elle est ouverte et ferm\'ee.
\item \label{lem-implicite3}  La surjection canonique $\ol{f_{\Top}}: \FX(A)_{\Top}\to \im(f_{\Top})$ admet assez de sections \emph{(cf.\ 2.1)}.
\end{numlist}
}

\smallskip

\dem 
\noindent\ref{lem-implicite1} Le cas \'etale (c'est-\`a-dire $d=0$) est une cons\'equence directe de la propri\'et\'e
hens\'elienne de $A$; voir \cite[Proposition 3.3.2]{GGMB}. Le cas g\'en\'eral s'en d\'eduit ais\'ement \`a l'aide d'une factorisation locale 
$$\FX \xrightarrow{\text{\;\'etale\;}}\FY\times\Aa^d\xrightarrow{\;\pr_{1}\;}\FY$$
 de $f$ au voisinage de $x$ (voir {\cite[I.4.4.2]{DG}}, ou \cite[17.11.4]{EGA4}).\smallskip

\noindent\ref{lem-implicite2} D'apr\`es \ref{lem-implicite1}, l'image de $f_{\Top}$ est r\'eunion de boules de rayon $J$; comme ces
boules forment une partition de $\FY(A)_{\Top}$ et sont ouvertes et ferm\'ees, \ref{lem-implicite2} en r\'esulte.\smallskip

\noindent\ref{lem-implicite3} est cons\'equence imm\'ediate de \ref{lem-implicite1} et \ref{lem-implicite2}.\qed

\subsection{Changement d'anneau}\label{SsecChgtAnn}
Soit $u:A\to K$ un homomorphisme continu (non n\'ecessairement local) d'anneaux locaux topologiques. (Le cas qui nous int\'eressera 
dans la suite est celui o\`u $A$ est un anneau de valuation et $K$ son corps des fractions.) 

Si $\FX$ est un $A$-sch\'ema localement de type fini, $\FX(K)$ s'identifie \`a $\FX_{K}(K)$ o\`u $\FX_{K}$ est le $K$-sch\'ema
d\'eduit de $\FX$ par changement de base. \`A ce titre, $\FX(K)$ est muni d'une topologie, et l'on a une application
 naturelle $\FX(A)\to \FX(K)$.
 
 \medskip
 
\noindent{\bf Lemme 3.4} 
{\it \begin{numlist}
\item \label{LemChgtAnn1} Avec les hypoth\`eses et notations de \rref{SsecChgtAnn}, l'application naturelle $\FX(A)\to \FX(K)$ est continue. 
\item \label{LemChgtAnn2} Si de plus $u$ est injectif et fait de $A$ un sous-anneau ouvert de $K$, alors cette application 
est ouverte; en particulier, c'est un plongement topologique ouvert si $\FX$ est s\'epar\'e.
\end{numlist}
}

\smallskip

\dem les  assertions \ref{LemChgtAnn1} et \ref{LemChgtAnn2} sont imm\'ediates 
par r\'eduction au cas affine.
\qed

\subsection{Valuations et valeurs absolues: rappels et conventions}\label{ssecVal}

Dans cet article, une \emph{valeur absolue} sur un corps $K$ est la donn\'ee d'un
groupe ab\'elien totalement 
ordonn\'e $(\Gamma,\cdot,1,\leq)$, not\'e multiplicativement et augment\'e 
d'un plus petit \'el\'e\-ment not\'e $0$, 
et d'une application
$$ \begin{array}{rrcl}
\abs: &K&\ffl&\Gamma\cup\{0\}\\
&z&\longmapsto&\abs(z)\quad\text{(\'egalement not\'e $\vert z\vert$)}
\end{array}$$
avec les propri\'et\'es habituelles: $\abs^{-1}(0)=\{0\}$, $\vert zz'\vert=\vert z\vert\,\vert z'\vert$, et:
\begin{itemize}
\item ou bien (cas ultram\'etrique) $\vert z+z'\vert\leq\max(\vert z\vert,\vert z'\vert)$, donc \og$\abs$\fg\ est une 
valuation (avec les conventions multiplicatives), d'anneau $O_{\abs}=\{z\in K\vert\; \abs(z)\leq1\}$;
\item ou bien (cas archim\'edien) $\Gamma{\subset}\RR_{>0}$ et \og$\abs$\fg\ est une valeur absolue ar\-chi\-m\'e\-dienne. Dans ce cas, nous conviendrons de poser $O_{\abs}:=K$.
\end{itemize}

Un \emph{corps valu\'e} est un corps muni d'une valeur absolue. Un corps valu\'e $(K,\Gamma,\abs)$ est un corps topologique: une
base de voisinages de $z_{0}\in K$ est donn\'ee par les boules $B(z_{0}, \rho):=\{z\in K\vert\; \abs(z-z_{0})<\rho\}$, o\`u $\rho$ 
parcourt $\Gamma$. Le compl\'et\'e de $K$ pour cette topologie est not\'e $\wh{K}$; c'est un corps valu\'e, qui dans le
cas archim\'edien est isomorphe \`a $\RR$ ou \`a $\CC$.

Un corps valu\'e $(K,\abs)$ est dit \emph{hens\'elien} si pour toute extension finie $L$ de $K$, il existe 
une \emph{unique} valeur 
absolue sur $L$ prolongeant celle de $K$. (Il existe toujours au moins un tel prolongement). Il est dit \emph{admissible} s'il 
est hens\'elien et si $\wh{K}$ est une extension s\'eparable de $K$. (Bien entendu, cette derni\`ere condition est satisfaite si $K$ est parfait, ou s'il est complet.)

Un corps valu\'e \emph{ultram\'etrique} $(K,\abs)$ est hens\'elien 
si et seulement si $O_{\abs}$ est un anneau local hens\'elien \cite[th. 4.1.3]{EP}. 

Un corps valu\'e  \emph{archim\'edien} $(K,\abs)$ est hens\'elien (et admissible!) si 
et seulement si $K$ est alg\'ebriquement clos (lorsque $\wh{K}\cong\CC$) ou r\'eel clos (lorsque $\wh{K}\cong\RR$): voir par exemple \cite[(26.8) et (26.9)]{E}. Dans les deux cas, il revient au m\^eme de dire que $K$ est al\-g\'e\-bri\-que\-ment ferm\'e dans $\wh{K}$.
\medskip

\noindent Sur un corps valu\'e hens\'elien, l'\'enonc\'e 3.3 admet la  variante suivante:

\medskip

\noindent{\bf Proposition 3.5} 
\textup{(Fonctions implicites, cas d'un 
corps valu\'e)} {\it Soit  $K$ un corps 
valu\'e hens\'elien et 
soit $f:X\to Y$ un morphisme \emph{lisse} de $K$-sch\'emas de type fini, de dimension
relative $d$. Soient $x\in X(K)$ et $y=f(x)$. 
\begin{numlist}
\item\label{lem-implicite-corps1} Il existe un voisinage $\Omega$ de $0$ dans $K$, et des voisinages $U\subset X(K)_{\Top}$
et $V\subset Y(K)_{\Top}$ de $x$ et $y$ 
respectivement, tels que $f_{\Top}(U)= V$ et que l'application $U\to V$ induite soit 
une fibration triviale de fibre $\Omega^d$.

En particulier, $f_{\Top}$ est ouverte (et est un hom\'eomorphisme local 
si $f$ est \'etale), et la surjection canonique 
$X(K)_{\Top}\to\im\,(f_{\Top})$ admet assez de sections locales.
\item\label{lem-implicite-corps2}
Si $K$ est ultram\'etrique, la surjection canonique 
$X(K)_{\Top}\to\im\,(f_{\Top})$ admet assez de sections.
\end{numlist}
}

\smallskip

\dem La partie \ref{lem-implicite-corps1} se d\'eduit facilement  de 3.3 dans le cas ultram\'etrique. Supposons $K$  archim\'edien. Si $K$ est complet, 
alors $K=\RR$ ou $\CC$ et il s'agit du th\'eor\`eme des fonctions implicites classique. 
En g\'en\'eral, il suffit de traiter le cas 
\'etale (m\^eme argument que dans 3.3\,\ref{lem-implicite1}). On consid\`ere alors le morphisme
$\wh{f}:X_{\wh{K}}\to Y_{\wh{K}}$ d\'eduit de $f$ par changement de base. Consid\'erons 
le diagramme commutatif d'espaces topologiques
$$\xymatrix{X(K)\ar[d]_{f_{\Top}}\ar@{}[r]|{\subset} & X(\wh{K})\ar[d]_{\wh{f}_{\Top}}\\
Y(K)\ar@{}[r]|{\subset} & Y(\wh{K})}$$
o\`u $\wh{f}_{\Top}$ est un hom\'eomorphisme local d'apr\`es le cas complet. 

Montrons que $X(K)=\wh{f}^{-1}_{\Top}(Y(K))$. 
Si $y\in Y(K)$, le $K$-sch\'ema $X_{y}=f^{-1}(y)$ est 
fini puisque $f$ est \'etale; 
en particulier ses corps r\'esiduels sont alg\'ebriques 
sur $K$, de sorte que $X_{y}(K)=X_{y}(\wh{K})$  puisque $K$ est 
al\-g\'e\-bri\-que\-ment ferm\'e dans $\wh{K}$, d'o\`u 
notre assertion.

\noindent\ref{lem-implicite-corps2} 
Supposons $K$ ultram\'etrique. En  vertu de la proposition 7.4, 
$\im(f_{\Top})$ est un espace ultraparacompact; comme l'application envisag\'ee a assez de sections locales, il s'ensuit facilement qu'elle a assez de sections.
\qed

\medskip

\noindent{\bf Proposition 3.6.} 
{\it Soit $K$ un corps valu\'e \emph{admissible}. 
\begin{numlist}
\item\label{prop-propre1} Soit $X$ un $K$-sch\'ema de type fini. Alors $X(K)$ est 
dense dans $X(\wh{K})_{\Top}$.
\item\label{prop-propre2} Soit $f:X\to Y$ un morphisme \emph{propre} de $K$-sch\'emas
de type fini. Alors l'image de $f_{\Top}:X(K)_{\Top}\to Y(K)_{\Top}$ est ferm\'ee.
\end{numlist} 
}

\smallskip

\dem 
dans le cas ultram\'etrique, voir \cite[1.2.1]{MB} pour \ref{prop-propre1}, et  \cite[1.3]{MB}  pour \ref{prop-propre2}.

Supposons $K$  archim\'edien, et montrons \ref{prop-propre1}. On peut supposer $X$ r\'eduit; il est alors r\'eunion 
finie de sous-$K$-sch\'emas localement ferm\'es \emph{lisses} ($K$ est de  caract\'eristique nulle). Ceci nous ram\`ene  
au cas o\`u $X$ est lisse, puis, quitte \`a le remplacer par un ouvert convenable, au cas o\`u l'on a un $K$-morphisme
\'etale $\pi: X\to \Aa^n_{K}$. Soit $\Omega\subset X(\wh{K})_{\Top}$ un ouvert non vide, et montrons que $\Omega\cap X(K)\neq\emptyset$. 
L'application induite $\wh{\pi}_{\Top}:X(\wh{K})\to \wh{K}^n$ est ouverte d'apr\`es 3.5, de sorte que 
$\wh{\pi}_{\Top}(\Omega)$ est ouvert dans $\wh{K}^n$ et contient donc un point $z\in K^n$ puisque $K$ est dense dans $\wh{K}$. 
On a donc $z=\wh{\pi}(y)$ pour un point $y\in \Omega$, mais on a n\'ecessairement $y\in X(K)$ par le m\^eme argument que dans 
la preuve de 3.5 ($\pi$ est quasi-fini, $\pi(y)$ est $K$-rationnel et $K$ est alg\'ebriquement ferm\'e 
dans $\wh{K}$): nous avons  bien trouv\'e un point de $X(K)\cap\Omega$.

Montrons enfin  \ref{prop-propre2} dans le cas archim\'edien; 
l'argument est encore ana\-logue \`a celui de 3.5. Comme $f$ est propre et $\wh{K}$ localement compact, 
$\wh{f}_{\Top}:X(\wh{K})\to Y(\wh{K})$ est topologiquement propre, de sorte que son image est fer\-m\'ee. Donc, si 
$y\in Y(K)$ est adh\'erent \`a $\im\,(f_{\Top})$, on a $y\in\im(\wh{f}_{\Top})$, d'o\`u $X_{y}(\wh{K})\neq\emptyset$. 
L'assertion \ref{prop-propre1} appliqu\'ee \`a $X_{y}$ implique alors que $X_{y}(K)\neq\emptyset$, et l'on conclut que $y\in\im(f_{\Top})$.\qed

\section{Produits et produits restreints}\label{SecProdRestr}

\subsection{Produits d'anneaux locaux topologiques}\label{SsecProduits}
On fixe un ensemble $V$ et une famille $(A_{v})_{v\in V}$, index\'ee par $V$, d'anneaux 
locaux topologiques (\ref{ssecAnnLocTop}). On pose
$A:=\prod_{v\in V}A_{v}$; c'est un anneau topologique, pour la topologie produit. 

Soit $\FX$ un $A$-sch\'ema. Pour chaque $v\in V$, on obtient par changement de base 
(par la projection $A\to A_{v}$) un 
$A_{v}$-sch\'ema not\'e $\FX_{v}$, et $\FX_{v}(A_{v})$ s'identifie \`a $\FX(A_{v})$. On 
a une application naturelle 
$$\FX(A)\ffl\prod_{v\in V}\FX_{v}(A_{v}). \eqno{(4.1)}
$$ 
Cette application est bijective si $\FX$ est 
quasi-compact et quasi-s\'epar\'e \cite[Theorem 3.6]{ConAdelic}%
\footnote{L'argument utilise le fait que les anneaux $K_{v}$ et $O_{v}$ ($v\in V$) sont 
tous locaux; cette restriction n'est 
en fait pas n\'ecessaire, cf. \cite[Theorem 1.3]{Bt} mais il s'agit d'un r\'esultat
bien plus profond.}%
, et notamment si $\FX$ est \emph{de pr\'esentation finie} sur $A$. Dans ce dernier cas, on peut donc munir $\FX(A)$ de la topologie 
produit des topologies naturelles sur les facteurs du membre de droite de (4.1). 
On obtient un espace topologique not\'e $\FX(A)_{\Top}$. 
Cette construction est fonctorielle en $\FX$, commute aux produits fibr\'es, transforme les immersions ferm\'ees en plongements topologiques ferm\'es et les 
immersions ouvertes en plongements topologiques (non ouverts en g\'en\'eral). En outre, l'espace $\FX(A)_{\Top}$ est s\'epar\'e si $\FX$ est s\'epar\'e.
Si $\FX=\Aa^n_{A}$, l'espace $\FX(A)_{\Top}$ s'identifie \`a $A^n$ avec sa topologie naturelle.

\subsection{Produits restreints de corps valu\'es}\label{ssecProdRestr} 

On fixe un ensemble d'indices $V$ et une famille  
$\soul{K}=(K_{v},\abs_{v})_{v\in V}$ de corps valu\'es. 
Pour $v\in V$ et $z\in K_{v}$, on utilisera aussi 
la notation $\vert z\vert_{v}$ pour $\abs_{v}(z)$.

\medskip

\noindent{\bf Remarque 4.1.} Dans les applications classiques, les $K_{v}$ sont les compl\'e\-t\'es d'un m\^eme corps $K$ pour une famille $V$ de valuations. 
Cependant, $K$ ne joue aucun r\^ole dans les g\'en\'eralit\'es sur la topologie ad\'elique, et  n'interviendra pas avant la section \ref{ssec-carac-nulle2}; d'autre part, l'hypoth\`ese que les $K_{v}$ soient complets peut le plus souvent \^etre remplac\'ee par une hypoth\`ese hens\'elienne. 

\medskip

On note $V_{\infty}$ l'ensemble des $v\in V$ tels que $K_{v}$ soit archim\'edien.
Pour $v\notin V_{\infty}$, $\abs_{v}$ est donc une valuation 
dont l'anneau sera not\'e $O_{v}$. Pour $v\in V_{\infty}$, on posera $O_{v}=K_{v}$.

Pour $S\subset V$ fini, on consid\`ere  l'anneau topologique produit
$$\ade_{\soul{K},S}:=\prod_{v\in S}K_{v}\times\prod_{v\notin S}O_{v}.$$

Pour $S\subset S'$ (tous deux finis), $\ade_{\soul{K},S}$ est un sous-anneau ouvert de $\ade_{\soul{K},S'}$ (et $\ade_{\soul{K},S'}$ est un localis\'e 
de $\ade_{\soul{K},S}$). Le \emph{produit restreint} des $K_{v}$, relativement aux $O_{v}$,  est par d\'efinition l'anneau topologique 
$$\ade_{\soul{K}}:=\varinjlim_{S\text{ fini }\subset V}\ade_{\soul{K},S}.$$
On peut le voir comme un sous-anneau de $\prod_{v\in V}K_{v}$ mais sa topologie 
n'est pas induite par 
la topologie produit. Lorsque l'on raisonne sur
la topologie de $\ade_{\soul{K}}$, il est utile de se souvenir que chaque anneau produit $\ade_{\soul{K},S}$ s'identifie \`a un sous-anneau ouvert de $\ade_{\soul{K}}$.

\subsection{Points des sch\'emas \`a valeurs dans un produit
restreint; topologie ad\'elique}\label{ssecPtsAd}
Dans tout le \S \ref{ssecPtsAd}, on fixe 
$\soul{K}=(K_{v},\abs_{v})_{v\in V}$ comme 
dans \ref{ssecProdRestr}. Si $X$ est un
$\ade_{\soul{K}}$-sch\'ema de pr\'esentation finie, on se propose de 
construire une topologie naturelle, dite \emph{ad\'elique}, 
sur l'ensemble
$X(\ade_{\soul{K}})$, selon une m\'ethode due \`a  Weil \cite{W} dans le cas des ad\`eles d'un corps global. Nous suivons la pr\'esentation de 
\cite{ConAdelic} auquel nous renvoyons pour les d\'etails. 
Noter qu'il n'est pas n\'ecessaire, m\^eme dans 
le cas classique des compl\'et\'es d'un corps 
global $K$, de se limiter aux  $\ade_{\soul{K}}$-sch\'emas 
provenant de $K$-sch\'emas de type fini.

\subsubsection{$S$-mod\`eles d'un sch\'ema; d\'efinition de la topologie ad\'elique. }\label{ssec-modeles} Soit $X$ 
un $\ade_{\soul{K}}$-sch\'ema de pr\'esentation finie. 
Comme  $\ade_{\soul{K}}$  est limite inductive filtrante des $\ade_{\soul{K},S}$, il existe d'apr\`es
 \cite[\S8]{EGA4} une partie finie $S$ de $V$, un $\ade_{\soul{K},S}$-sch\'ema $\FX_{S}$ de pr\'esentation finie et 
 un $\ade_{\soul{K}}$-isomorphisme $\ade_{\soul{K}}\otimes_{\ade_{\soul{K},S}}\FX_{S}\cong X$; on dira que $\FX_{S}$ est un \emph{$S$-mod\`ele} de $X$.
 Dans ces conditions, l'application naturelle 
$$
\varinjlim_{S'\supset S} \FX_{S}(\ade_{\soul{K},S'}) \ffl X(\ade_{\soul{K}}) \eqno{(4.2)}
$$ 
est bijective (o\`u $S'$ parcourt les parties finies de $V$ contenant $S$). On munit 
alors $X(\ade_{\soul{K}})$ de 
la \emph{topologie limite inductive} 
d\'eduite de (4.2), cf. \cite[App. II.I]{D}.
On notera $X(\ade_{\soul{K}})_{\Top}$ l'espace topologique 
ainsi obtenu; il est ind\'ependant du choix 
du mod\`ele $\FX_{S}\to\Spec(\ade_{\soul{K},S})$.
Si $Y$ d\'esigne un espace topologique, 
une application $X(\ade_{\soul{K}})_{\Top} \to Y$ 
est continue si chaque application induite $X(\ade_{\soul{K}, S})_{\Top}
\to Y$ est continue.

\subsubsection{Propri\'et\'es de la topologie ad\'elique. }\label{ssecPropTopAd}
Gardons les notations de \ref{ssecPtsAd}. Le syst\`eme inductif des espaces $\FX_{S}\left(\ade_{\soul{K},S'}\right)_{\Top}$, dont 
$X(\ade_{\soul{K}})_{\Top}$ est la colimite, admet la propri\'et\'e agr\'eable suivante, cons\'equence de 3.4: 
pour $S\subset S'\subset S''$, l'application de transition
$$\FX_{S}\left(\ade_{\soul{K},S'}\right)_{\Top}\to \FX_{S}\left(\ade_{\soul{K},S''}\right)_{\Top}$$
 est continue et \emph{ouverte}, et est un plongement topologique ouvert si $\FX_{S}$ est s\'epar\'e (ce que l'on 
 peut toujours supposer lorsque $X$ est s\'epar\'e, 
 cf. \cite[8.10.4]{EGA4}). 
 Par suite, on a les m\^emes propri\'et\'es pour les applications naturelles $\FX_{S}\left(\ade_{\soul{K},S'}\right)_{\Top}\to X\left(\ade_{\soul{K}}\right)_{\Top}$. 
 L'\'etude d'un espace $X\left(\ade_{\soul{K}}\right)_{\Top}$ se fait donc g\'en\'eralement en deux \'etapes:
 \begin{itemize}
\item choix d'un mod\`ele $\FX_{S}$ adapt\'e, et \'etude des espaces  $\FX_{S}\left(\ade_{\soul{K},S'}\right)_{\Top}$: ceux-ci sont des espaces produits,
justiciables d'\'enonc\'es g\'en\'eraux tels que ceux de  \ref{SsecTopProduit};
\item passage \`a la limite, utilisant le fait que les $\FX_{S}\left(\ade_{\soul{K},S'}\right)_{\Top}$ sont ouverts dans $X\left(\ade_{\soul{K}}\right)_{\Top}$.
\end{itemize}

La construction $X\mapsto X(\ade_{\soul{K}})_{\Top}$ est fonctorielle en $X$, commute aux pro\-duits fibr\'es, transforme les immersions ferm\'ees en
plongements topologiques ferm\'es, les immersions ouvertes en plongements topologiques, les sch\'emas s\'epar\'es en espaces s\'epar\'es; 
en outre, $(\Aa^n_{\ade_{\soul{K}}})_{\Top}$ s'identifie \`a $(\ade_{\soul{K}})^n$, pour tout $n\in\NN$.

\subsubsection{Groupes et torseurs. }\label{ssec:groupes-torseurs}\strut
La compatibilit\'e aux produits implique notamment que si $G$ est un  $\ade_{\soul{K}}$- sch\'ema en groupes de pr\'esentation finie op\'erant sur $X$,
alors $G(\ade_{\soul{K}})_{\Top}$ est un groupe topologique op\'erant contin\^ument sur $X(\ade_{\soul{K}})_{\Top}$.

L'objet de cet article est d'\'etudier la question suivante,
en analogie avec le cas d'un corps valu\'e \cite{GGMB}.

\bigskip

\noindent{\bf Question.} Soient $G$ un $\ade_{\soul{K}}$-sch\'ema en groupes de pr\'esentation finie et  ${f: X \to Y}$ un $G$-torseur,  o\`u 
$X$ et $Y$ sont des $\ade_{\soul{K}}$-sch\'emas de pr\'esentation finie. Consid\'erons l'application
$f_{\Top}: X(\ade_{\soul{K}})_{\Top}\to Y(\ade_{\soul{K}})_{\Top}$ induite sur les points ad\'eliques, et son image $I\subset Y(\ade_{\soul{K}})_{\Top}$. 
\begin{itemize}
\item \`A quelles conditions (sur $G$, notamment) $I$ est-elle localement ferm\'ee (resp. ouverte, ferm\'ee) dans $Y(\ade_{\soul{K}})_{\Top}$?
\item L'application induite $X(\ade_{\soul{K}})_\Top \to I$ est-elle un $G_{\Top}$-fibr\'e principal? 
\end{itemize}
\subsubsection{$S$-mod\`eles d'un morphisme. }\label{ssec-modeles-morphisme} 
Pour simplifier, nous ne con\-si\-d\'e\-re\-rons 
d\'esormais que des $\ade_{\soul{K}}$-sch\'emas \emph{s\'epar\'es}, 
et des mod\`eles s\'epar\'es de ceux-ci. 
\smallskip

Consid\'erons un morphisme $f:X\to Y$ de $\ade_{\soul{K}}$-sch\'emas s\'epar\'es de 
pr\'e\-sen\-ta\-tion finie\footnote{Rappelons qu'un tel morphisme est 
s\'epar\'e  \cite[25.21.13, Tag 01KV]{St}, \cite[(5.3.1.(v))]{EGA1}.
}.
Si $S$ est une partie finie de $V$, nous appellerons \emph{$S$-mod\`ele} de $f$ 
la donn\'ee de $S$-mod\`eles $\FX_{S}$ et $\FY_{S}$ de $X$ et $Y$ et d'un $\ade_{\soul{K},S}$-morphisme $F_{S}:\FX_{S}\to\FY_{S}$ induisant $f$ par 
le changement de base $\ade_{\soul{K},S}\to\ade_{\soul{K}}$. D'apr\`es  \cite[\S8]{EGA4}, un tel mod\`ele existe toujours (pour $S$ convenable); en outre, 
si $F^{(j)}_{S_{j}}:\FX_{S_{j}}\to\FY_{S_{j}}$  ($j=1,2$) sont deux mod\`eles de $f$, il existe $S$ fini contenant $S_{1}$ et $S_{2}$ et des
isomorphismes $\FX_{S_{1}}\otimes_{\ade_{\soul{K},S_{1}}}\ade_{\soul{K},S} \flis\FX_{S_{2}}\otimes_{\ade_{\soul{K},S_{2}}}\ade_{\soul{K},S}$
et $\FY_{S_{1}}\otimes_{\ade_{\soul{K},S_{1}}}\ade_{\soul{K},S} \flis\FY_{S_{2}}\otimes_{\ade_{\soul{K},S_{2}}}\ade_{\soul{K},S}$ compatibles, au sens \'evident,
avec $F^{(1)}_{S_{1}}$ et $F^{(2)}_{S_{2}}$.
\subsection{\hskip-0.5em Le diagramme des images }\label{sec-diag-images}%

Dans ce num\'ero, on  utilise les conventions suivantes
pour les applications continues: 
\begin{gather}
 \notag \xymatrixcolsep{.5cm}
\xymatrix{{}\ar@{{}->>}[r]&{}}  \raisebox{.5ex}{= surjection}\qquad
\xymatrix{{}\ar@{^{(}->}[r]&{}} \raisebox{.5ex}{= injection}\qquad
\xymatrix{{}\ar^{\cong}[r]&{}} \raisebox{.5ex}{= hom\'eomorphisme}
\\ \notag \xymatrixcolsep{.5cm}
\xymatrix{{}\ar@{^{(}->}|(.4){\bullet}[r]&{}}  \raisebox{.5ex}{= plongement topologique} \quad\qquad
\xymatrix{{}\ar@{^{(}->}|(.4){\circ}[r]&{}} \raisebox{.5ex}{= plongement topologique ouvert.}
\end{gather}

Soit $(K_{v})_{v\in V}$ comme dans \rref{ssecProdRestr}, et soit  $f: X \to Y$ un 
morphisme de $\ade_{\soul{K}}$-sch\'emas s\'epar\'es
de pr\'esentation finie. Soient $S\subset V$ fini et $F_{S}:\FX_{S}\to\FY_{S}$  un $S$-mod\`ele (s\'epar\'e) de $f$. 
Pour abr\'eger, posons, pour tout $v\in V$,
$$O'_{v}:=\begin{cases} K_{v} &\text{si $v\in S$}\\ O_{v}& \text{si $v\notin S$}\end{cases}$$
de sorte que $\ade_{\soul{K},S}=\prod_{v\in V}O'_{v}$; en outre,  le $O'_{v}$-morphisme  $F_{S,v}:\FX_{S,v}\to\FY_{S,v}$ 
d\'eduit de $F_{S}$ 
sera not\'e plus simplement  $F_{v}:\FX_{v}\to\FY_{v}$, et $F_{v,\Top}$ d\'esignera l'application continue induite
$\FX_{S}(O'_{v})_{\Top} \to \FY_{S}(O'_{v})_{\Top}$. 

Notre but est d'\'etudier l'image de $f_{\Top}:X(\ade_{\soul{K}})_{\Top}\to Y(\ade_{\soul{K}})_{\Top}$ et ses liens avec les images des
applications locales (i.e. $f_{v,\Top}$ et $F_{v,\Top}$) et $S$-ad\'eliques (i.e. $F_{S,\Top}$). La situation est r\'esum\'ee
par le diagramme commutatif d'espaces topologiques
\begin{gather}  \nonumber 
\raisebox{3cm}
{\xymatrixcolsep{1.5cm}
$\xymatrix{%
\prod\limits_{v\in V}\FX_{S}(O'_{v})_{\Top}
\ar@{{}->>}^{\prod_{v}\ol{F}_{v,\Top}}[r] &
\prod\limits_{v\in V}\im\,(F_{v,\Top})\;
\ar@{^{(}->}|{\bullet}[r] \ar@{}|{\txt{\small(A)}}[dr] &
\prod\limits_{v\in V}\FY_{S}(O'_{v})_{\Top}\\
{\FX_{S}(\ade_{\soul{K},S})_{\Top}}_{\strut} \;
\ar_(.4){\cong}^(.4){\varphi_{X}}[u] \ar@{^{(}->}|(.6){\circ}_(.6){j_{X}}[d]+<0ex,2.5ex> \ar@{{}->>}^{\ol{F}_{S,\Top}}[r]&
\im\,(F_{S,\Top})_{\strut}\;
\ar_(.4){\cong}^(.4){\varphi_{\mathrm{Im}}}[u] \ar@{^{(}->}|(.6){\bullet}_(.6){j_{\mathrm{Im}}}[d]+<0ex,2.5ex> \ar@{}|{\txt{\small(B)}}[dr] 
\ar@{^{(}->}|(.46){\bullet}^{u_{F_{S}}}[r] &
{\FY_{S}(\ade_{\soul{K},S})_{\Top}}_{\strut} \;
\ar_(.4){\cong}^(.4){\varphi_{Y}}[u] \ar@{^{(}->}|(.6){\circ}_(.6){j_{Y}}[d]+<0ex,2.5ex> \\
 {X(\ade_{\soul{K}})_{\Top}}_{\strut}^{\strut} \;
\ar@{{}->>}^{\ol{f}_{\Top}}[r]\ar@{^{(}->}_{i_{X}}^{\strut}[d] & 
 \im\,(f_{\Top})_{\strut}^{\strut}\;
\ar@{^{(}->}_{i_{\mathrm{Im}}}[d] \ar@{}|{\txt{\small(C)}}[dr]
\ar@{^{(}->}|(.46){\bullet}^{u_{f}}[r] &
 {Y(\ade_{\soul{K}})_{\Top}}_{\strut}^{\strut} \;
\ar@{^{(}->}_{i_{Y}}[d] & \qquad (4.3) \\
\prod\limits_{v\in V} X_{v}(K_{v})_{\Top}
\ar@{{}->>}^{\prod_{v}\ol{f}_{v,\Top}}[r] &
\prod\limits_{v\in V}\im\,(f_{v,\Top}) \;
\ar@{^{(}->}|{\bullet}[r] &
 \prod\limits_{v\in V} Y_{v}(K_{v})_{\Top}
\save "3,1"+<-7ex,0ex>."3,3"*\frm<8pt>{.}\restore}$
}
 \end{gather}
dans lequel:
\begin{itemize}
\item chaque image (colonne du milieu) est munie de la topologie induite par l'espace d'arriv\'ee (\`a sa droite). 
Les quatre fl\`eches horizontales du type $\xymatrix{{}\ar@{^{(}->}|{\bullet}[r]&{}}$ sont donc des
\emph{plongements topologiques} (cf. 2.2.\,\ref{lemTopProduit1});
\item les fl\`eches $\varphi_{X}$ et $\varphi_{Y}$ sont des hom\'eomorphismes par d\'efinition des 
topologies sur $\FX_{S}(\ade_{\soul{K},S})$ et $\FY_{S}(\ade_{\soul{K},S})$ (cf. \ref{SsecProduits}); 
\item $\varphi_{\mathrm{Im}}$ est \'egalement un ho\-m\'eo\-mor\-phisme d'apr\`es 2.3\,\ref{corTopProduit1}, et 
le carr\'e (A) est (trivialement) cart\'esien;
\item les deux fl\`eches descendantes du type $\xymatrix{{}\ar@{^{(}->}|{\circ}[r]&{}}$ sont des \emph{plongements topologiques ouverts} (\ref{ssecPropTopAd});
\item puisque $j_{Y}$, $u_{F_{S}}$ et $u_{f}$ sont des plongements topologiques, $j_{\mathrm{Im}}$ en est un;
\item le sous-diagramme form\'e des deux derni\`eres lignes ne d\'epend pas de $S$ ni du mod\`ele $F_{S}$; 
\item en g\'en\'eral, $i_{X}$, $i_{\mathrm{Im}}$ et $i_{Y}$ ne sont pas  des plongements topologiques.
\end{itemize} \medskip

Lorsque $S$ varie, chacun des ensembles de la troisi\`eme ligne 
est r\'eunion filtrante des images des ensembles du dessus, en sorte
que les deuxi\`eme et troisi\`eme lignes donnent naissance au diagramme

\begin{gather} \nonumber 
\raisebox{1cm}
{\xymatrixcolsep{1.5cm}
$\xymatrix{%
\varinjlim\limits_{S}\:{\FX_{S}(\ade_{\soul{K},S})_{\Top}}_{\strut} \;
\ar^(.65){\cong}_(.65){\psi_{X}}[d]+<0ex,2.5ex> \ar@{{}->>}[r]&
\varinjlim\limits_{S}\:\im\,(F_{S,\Top})_{\strut}\;
\ar@{^{(}->>}_(.65){\psi_{\mathrm{Im}}}[d]+<0ex,2.5ex> 
\ar@{^{(}->}[r] &
\varinjlim\limits_{S}\:{\FY_{S}(\ade_{\soul{K},S})_{\Top}}_{\strut} \;
\ar^(.65){\cong}_(.65){\psi_{Y}}[d]+<0ex,2.5ex> \\
 {X(\ade_{\soul{K}})_{\Top}}_{\strut}^{\strut} \;
\ar@{{}->>}^{\ol{f}_{\Top}}[r] & 
 \im\,(f_{\Top})_{\strut}^{\strut}\;
\ar@{^{(}->}|(.46){\bullet}^{u_{f}}[r] &
 {Y(\ade_{\soul{K}})_{\Top}}_{\strut}^{\strut} \;
 & (4.4)
\save "2,1"+<-7ex,0ex>."2,3"*\frm<8pt>{.}\restore}$
}
 \end{gather}
o\`u les fl\`eches $\psi_{X}$ et $\psi_{Y}$ sont des hom\'eomorphismes
par d\'efinition des topologies ad\'eliques, alors que $\psi_{\mathrm{Im}}$ est, 
a priori, seulement une bijection continue.

La difficult\'e, dans l'\'etude de l'application ad\'elique $f_{\Top}$ (troisi\`eme ligne), r\'eside en g\'en\'eral dans le passage \`a la 
limite \`a partir des propri\'et\'es de $F_{S,\Top}$ (deuxi\`eme ligne). Ces derni\`eres  se ram\`enent \`a celles 
de la premi\`ere ligne: la proposition ci-dessous traite le cas o\`u $F_{S}$ est lisse.

\medskip

\noindent{\bf Proposition 4.2.} 
{\it Avec les hypoth\`eses et notations de \rref{sec-diag-images}, on suppose que $V_{\infty}$ est \emph{fini}, 
que $F_{S}:\FX_{S}\to\FY_{S}$ est \emph{lisse} et que, pour tout $v\in V$, le corps valu\'e $K_{v}$ est \emph{hens\'elien}. 
\begin{romlist}
\item\label{prop-lisse-S-adelique1} La surjection 
$\ol{F}_{S,\Top}:
\FX_{S}(\ade_{\soul{K},S})_{\Top}\to\im\,(F_{S,\Top})$ 
a assez de sections locales, et assez de sections si 
$V_\infty=\emptyset$.
\item\label{prop-lisse-S-adelique2} Le plongement $u_{F_{S}}:\im\,(F_{S,\Top})\inj \FY_{S}(\ade_{\soul{K},S})_{\Top}$ est \emph{localement ferm\'e}.
\item\label{prop-lisse-S-adelique3}  Supposons de plus que, pour presque tout $v\in V$, l'application $F_{v,\Top}:\FX_{S}(O'_{v})_{\Top}\to\FY_{S}(O'_{v})_{\Top}$ 
soit \emph{surjective}. Alors le plongement $u_{F_{S}}$ est ouvert.
\end{romlist}
}

\smallskip

\dem il suffit d'\'etablir les propri\'et\'es correspondantes pour les fl\`eches de la 
ligne sup\'erieure du diagramme (4.3). 
On peut de plus supposer que $V_{\infty}\subset S$ (remarquer que $\ade_{\soul{K},S\cup V_{\infty}}=\ade_{\soul{K},S}$). 

Dans ces conditions, pour $v\in S$ (et donc $O'_{v}=K_{v}$), la proposition 3.5 nous dit que $\im\,(F_{v,\Top})$ est
un ouvert de $\FY_{S}(O'_{v})_{\Top}$ et que $\ol{F}_{v,\Top}$ admet assez de sections locales; si de 
plus $v\in S\smallsetminus V_{\infty}$, elle a m\^eme assez de sections d'apr\`es 3.5.\,\ref{lem-implicite-corps2}. D'autre part, 
pour $v\notin S$, on a $O'_{v}=O_{v}$, et 
la proposition 3.3 indique que $\im\,(F_{v,\Top})$ est ouvert et ferm\'e dans $\FY_{S}(O'_{v})_{\Top}$ et
que $\ol{F}_{v,\Top}$ admet assez de sections.
En vertu du corollaire 2.3\,\ref{corTopProduit3},
$\ol{F}_{S,\Top}:
\FX_{S}(\ade_{\soul{K},S})_{\Top}\to\im\,(F_{S,\Top})$ a assez de sections locales, d'o\`u 
la premi\`ere partie de l'assertion
\ref{prop-lisse-S-adelique1}; si de plus  $V_\infty= \emptyset$,  
le raffinement  2.3\,\ref{corTopProduit2} indique
que cette application a assez de sections. 
En outre, le lemme 2.2.\,\ref{lemTopProduit4}
montre que $\im\,(F_{S,\Top})$ est localement ferm\'e dans 
$\prod\limits_{v\in V}\FY_{S}(O'_{v})_{\Top}$, c'est-\`a-dire l'assertion 
\ref{prop-lisse-S-adelique2}.

Sous l'hypoth\`ese suppl\'ementaire de   \ref{prop-lisse-S-adelique3}, 
on a $\im\,(F_{v,\Top})=\FY_{S}(O'_{v})_{\Top}$ pour presque tout $v$, de sorte que $\prod_{v\in V}\im\,(F_{v,\Top})$
est bien ouvert d'apr\`es 2.2.\,\ref{lemTopProduit5}.\qed

\section{La propri\'et\'e \og presque  propre \fg }\label{secPP}

\subsection{Version locale }\label{ssecPPL}Soit $(K,v)$ un corps valu\'e; on pose $O=O_{v}$ et l'on consid\`ere un morphisme $f:\FX\to\FY$ de $O$-sch\'emas s\'epar\'es. 
Nous noterons PPL (\og presque-propret\'e locale\fg) la condition suivante sur $f$:
\medskip

\noindent $\mathrm{PPL}(O,f)$:  \qquad$\mathrm{Im}\bigl( \FX(O) \to \FY(O) \bigr) \;=\;
\mathrm{Im}\bigl( \FX(K) \to \FY(K) \bigr) \; \cap \;\FY(O)$
\medskip

\noindent o\`u l'on identifie $\FX(O)$ et $\FY(O)$ \`a leurs images respectives dans $\FX(K)$ et $\FY(K)$.

\medskip

\noindent{\bf Remarques 5.1.} 
\begin{numlist}
\item\label{remPPL1}  La condition (PPL) est trivialement v\'erifi\'ee si $v$ 
est archi\-m\'e\-dienne: dans ce cas, on a $O=K$ par convention. 
\item\label{remPPL1,5}  Un autre cas trivial, un peu plus int\'eressant, est
celui o\`u l'application $\FX(O) \to \FY(O)$ est surjective.
\item\label{remPPL2}  (PPL) est aussi satisfaite si $f$ est propre; plus g\'en\'eralement, 
si $f:\FX\to \FY$ v\'erifie (PPL) et si $g:\FY\to\FZ$ est propre, 
alors $g\circ f:\FX\to\FZ$ v\'erifie (PPL).
C'est l\`a une cons\'equence facile du crit\`ere valuatif de propret\'e.
\end{numlist}

\subsection{Version ad\'elique }\label{ssecPPA} 
Soit $(K_{v})_{v\in V}$ comme dans \ref{ssecProdRestr}, dont nous reprenons les notations. 
Consid\'erons un morphisme $f:X\to Y$ de $\ade_{\soul{K}}$-sch\'emas s\'epar\'es de pr\'e\-sen\-ta\-tion finie. 

\medskip

\noindent{\bf D\'efinition 5.2.} 
{\it Avec les notations ci-dessus, nous dirons 
que $f: X \to Y$ est \emph{presque propre} s'il v\'erifie la condition suivante:\medskip

\noindent $\mathrm{PP}(\ade_{\soul{K}},f)$: \; il
 existe une partie finie $S$ de $V$ et un $S$-mod\`ele $F_{S}:\FX_{S}\to\FY_{S}$ de $f$ tel que  pour presque
 tout $v\notin S$ la condition $\mathrm{PPL}(O_{v},F_{S,v})$ de \rref{ssecPPL} soit satisfaite, 
 o\`u $F_{S,v}: \FX_{S,v}\to\FY_{S,v}$  est le $O_{v}$-morphisme 
 induit par $F_{S}$.
 \medskip

\noindent Nous dirons en outre que $F_{S}:\FX_{S}\to\FY_{S}$ est un \emph{bon $S$-mod\`ele} de $f$ 
si la condition $\mathrm{PPL}(O_{v},F_{S,v})$ est v\'erifi\'ee \emph{pour tout} $v\notin S$.}

\medskip

\noindent{\bf Remarques 5.3.} 
\begin{numlist}
\item\label{rem_propre1} Si $f$ est presque propre, la condition \'enonc\'ee est en fait valable pour 
tout $S$ et tout $S$-mod\`ele de $f$. De plus, \'etant
donn\'e un  tel mod\`ele $F_{S}$, soit $S'$ la r\'eunion de $S$ et de l'ensemble (fini) 
des $v\notin S$ tels que  $\mathrm{PPL}(O_{v},F_{v})$ ne
soit pas satisfaite. Alors  le $S'$-mod\`ele $F_{S'}:\FX_{S'}\to\FY_{S'}$ d\'eduit 
de $F_{S}$ par changement de base est un bon $S'$-mod\`ele.
\item\label{rem_propre2} Si $f$ est propre, alors $f$ satisfait la propri\'et\'e (PP). En effet
si $F_{S}:\FX_{S}\to\FY_{S}$ est un $S$-mod\`ele de $f$, il existe $S'$ fini
contenant $S$ tel que le $S'$-mod\`ele d\'eduit de $F$ par changement de base soit
propre \cite[(8.10.5)(xii)]{EGA4}; en particulier, pour tout $v\notin S'$,
le $O_{v}$-morphisme $F_{v}:\FX_{v}\to\FY_{v}$ est propre et v\'erifie donc (PPL).
\item\label{rem_propre3} On a plus g\'en\'eralement l'analogue de 5.1\,\ref{remPPL2}: 
si $f:X\to Y$ v\'erifie (PP) et 
si $g:Y\to Z$ est propre, alors $g\circ f:X\to Z$ v\'erifie (PP). 
\end{numlist}

\medskip

Nous verrons plus loin des exemples de cette
propri\'et\'e, cf. (\ref{ssecLangWeil}).

\medskip

\noindent{\bf Proposition 5.4.} 
 {\it Avec les notations de \rref{ssecPPA}, on suppose que $f: X\to Y$ est 
presque propre, et l'on fixe un bon
$S$-mod\`ele  $F_{S}:\FX_{S}\to\FY_{S}$ de $f$, pour $S\subset V$ fini convenable. On  consid\`ere 
le diagramme des images (4.3) associé (voir Section 4.4)
 Alors:
\begin{romlist}
\item\label{lemme-pp-cart1} les carr\'es  \textup{(B)}  et \textup{(C)} du diagramme sont 
\emph{cart\'esiens};
\item\label{lemme-pp-cart2}  $j_{\im}: \im\,(F_{S,\Top})\to \im\,(f_{\Top})$ est un 
plongement topologique \emph{ouvert};
\item\label{lemme-pp-cart3} la bijection continue 
$${\psi_{\mathrm{Im}}}:\varinjlim\limits_{S'\supset S}\:\im\,(F_{S',\Top})\to  \im\,(f_{\Top})$$ 
du diagramme \textup{\ref{sec-diag-images}\,(4.4)} est un hom\'eomorphisme.
\end{romlist}
}

\smallskip

\dem \ref{lemme-pp-cart1} notons d'abord que puisque les fl\`eches horizontales sont 
des plongements topologiques, il revient au m\^eme de dire que (B) (resp. (C)) est
ensemblistement cart\'esien ou qu'il l'est topologiquement. 

Montrons d'abord que le rectangle $\left(\begin{smallmatrix}\text{\tiny B}\\\text{\tiny C}\end{smallmatrix}\right)$ 
est cart\'esien. Composant avec les isomorphismes 
en haut du diagramme (4.3), il suffit de voir que le diagramme (5.1) de gauche 

$$ 
\raisebox{1cm}{%
$\xymatrix{%
\prod\limits_{v\in V}\im\,(F_{v,\Top})\;
\ar@{^{(}->}|{\bullet}[r] \ar@{^{(}->}[d] 
&
\prod\limits_{v\in V}\FY_{S}(O'_{v})_{\Top} \ar@{^{(}->}[d] &
\im\,(F_{v,\Top})_{\strut}\;
\ar@{^{(}->}|{\bullet}[r] \ar@{^{(}->}[d] 
&
{\FY_{S}(O'_{v})_{\Top}}_{\strut} \ar@{^{(}->}[d] 
 \\
\prod\limits_{v\in V}\im\,(f_{v,\Top})^{\strut} \;
\ar@{^{(}->}|{\bullet}[r] &
 {\prod\limits_{v\in V} Y_{v}(K_{v})_{\Top}}^{\strut} &
 \im\,(f_{v,\Top})^{\strut} \;
\ar@{^{(}->}|{\bullet}[r] &
{Y_{v}(K_{v})_{\Top}}^{\strut} & (5.1)
}$}
$$
est cart\'esien; or celui-ci est le produit, index\'e par $v\in V$, des diagrammes de droite, qui sont cart\'esiens par d\'efinition d'un bon mod\`ele.

Montrons maintenant que (C) est cart\'esien si $f$ est presque propre. 
L'injection $u_{f}: \im\,(f_{\Top})\inj Y(\ade_{\soul{K}})$
est (au moins ensemblistement) limite inductive filtrante des 
$u_{F_{S}}: \im\,(F_{S,\Top})\inj \FY_{S}(\ade_{\soul{K}},S)$ lorsque 
$F_{S}$ parcourt les bons mod\`eles de $f$. L'assertion est donc cons\'equence 
de la pr\'ec\'edente par passage \`a la limite.

Enfin, puisque le rectangle $\left(\begin{smallmatrix}\text{\tiny B}\\\text{\tiny C}\end{smallmatrix}\right)$  et
le carr\'e (C) sont cart\'esiens, le carr\'e (B) l'est aussi.

Cette derni\`ere propri\'et\'e implique aussi l'assertion \ref{lemme-pp-cart2} puisque
$j_{Y}$ est un plongement topologique ouvert. On en d\'eduit enfin 
 \ref{lemme-pp-cart3}: en effet, $\im\,(f_{\Top})$ est r\'eunion filtrante des images 
 des $\im\,(F_{S',\Top})$, lesquelles sont ouvertes, de sorte que cette r\'eunion est 
 aussi une colimite topologique.
 \qed

\medskip

\noindent{\bf Proposition 5.5.} 
{\it Soient $(K_{v})_{v\in V}$ et $f: X \to Y$ comme dans \rref{ssecPPA}. 
On suppose que $f: X \to Y$ est \emph{presque propre (5.2)}.
On note $I$ l'image de $f_\Top: X(\ade_{\soul{K}})_\Top \to Y(\ade_{\soul{K}})_\Top$.

\begin{numlist} 
\item\label{prop_pp2} On suppose que pour tout $v\in V$ l'image de $f_v: X(K_v)_{\Top} \to Y(K_v)_{\Top}$ 
est ferm\'ee. Alors $I$ est   ferm\'ee.
\item\label{prop_pp1} On suppose que $V_{\infty}$ est fini, que $f$ est  \emph{lisse}, et que 
tous les corps valu\'es $K_{v}$ sont \emph{hens\'eliens}. Alors:
\begin{romlist}
\item\label{prop_pp1_1} La surjection $\ol{f}_{\Top}:X(\ade_{\soul{K}})_{\Top}\to I$ 
a assez de sections locales. 
En particulier, $f_\Top$ est ouverte sur son image (et donc stricte). 
Si $V_{\infty}=\emptyset$, alors $\ol{f}_{\Top}$ a assez de sections.
\item\label{prop_pp1_2} $I$ est \emph{localement ferm\'ee} dans $Y(\ade_{\soul{K}})_{\Top}$.

\item\label{prop_pp1_3} Soit $F_{S}:\FX_{S}\to\FY_{S}$ un $S$-mod\`ele de $f$, et
supposons de plus que, pour presque tout $v\in V$, 
l'application $F_{v,\Top}:\FX_{S}(O'_{v})_{\Top}\to\FY_{S}(O'_{v})_{\Top}$ soit
\emph{surjective} \emph{(condition qui ne d\'epend pas du choix du mod\`ele)}. 
Alors $I$ est \emph{ouverte} dans $Y(\ade_{\soul{K}})_{\Top}$.
\end{romlist}
\end{numlist}
}

\smallskip

\dem  \ref{prop_pp2} L'hypoth\`ese sur les images entra\^{\i}ne que, dans le diagramme \ref{sec-diag-images}\,(4.3), le 
plongement en bas \`a droite est ferm\'e. D'autre part, comme $f$ est presque propre, le carr\'e (C) est 
cart\'esien (proposition 5.4), d'o\`u la conclusion.\medskip

\noindent\ref{prop_pp1} 
Il existe $S\subset V$ fini et un bon $S$-mod\`ele \emph{lisse} $F_{S}: \FX_{S}\to \FY_{S}$ de $f$. Pour tout $S'$ fini contenant $S$, on 
notera $F_{S'}:\FX_{S'}\to \FY_{S'}$ le bon $S'$-mod\`ele d\'eduit de $F_{S}$ par changement de base, et 
$I_{S'}\subset \FY_{S'}(\ade_{\soul{K},S'})_{\Top}$ l'image de $F_{S'}$. 

Dans ces conditions, $Y(\ade_{\soul{K}})_{\Top}$ est la r\'eunion filtrante, index\'ee par $S'$, des ouverts  $\FY_{S'}(\ade_{\soul{K},S'})_{\Top}$. 
De plus, comme chaque $F_{S'}$ est un bon mod\`ele, la proposition 5.4 nous dit que $I_{S'}=I\cap\FY_{S'}(\ade_{\soul{K},S'})_{\Top}$. 
Les trois assertions r\'esultent donc des assertions correspondantes de la proposition 4.2 appliqu\'ee 
\`a chaque $F_{S'}$.\qed
\medskip

\noindent{\bf Corollaire 5.6.} 
{\it  Soient $K$, $(K_{v})_{v\in V}$ et $f: X \to Y$ comme dans \rref{ssecPPA}. 
On note $I$ l'image de $f_\Top: X(\ade_{\soul{K}})_\Top \to Y(\ade_{\soul{K}})_\Top$, et l'on suppose que:
\begin{itemize}
\item $f$ est \emph{propre et lisse};
\item $V_{\infty}$ est \emph{fini};
\item tous les corps valu\'es $K_{v}$ sont \emph{admissibles} \textup{(\ref{ssecVal})}.
\end{itemize}
Alors $I$ est ferm\'ee dans $Y(\ade_{\soul{K}})_\Top$, et  
l'application  induite $X(\ade_{\soul{K}})_\Top  \to I$ 
admet assez de sections locales.
}

\smallskip

\dem Pour tout $v\in V$, l'image de $f_{v,\Top}:X(K_{v})_\Top \to Y(K_{v})_\Top$ est ferm\'ee d'apr\`es 3.6.\,\ref{prop-propre2}, 
puisque $f$ est propre et $K_{v}$ admissible. Par 
suite $I$ est ferm\'ee (5.5.\,\ref{prop_pp2}). 
L'existence d'assez de sections locales 
r\'esulte de 5.5.\,\ref{prop_pp1}.\qed
\medskip

\subsection{Morphismes presque propres: le crit\`ere d'Oesterl\'e}\label{ssecLangWeil}

Un cas int\'eressant o\`u (PP) est v\'erifi\'ee a \'et\'e remarqu\'e
par Oesterl\'e dans le cas affine \cite[I.3.6]{O} et raffin\'e
par Conrad \cite[Th. 4.5]{ConAdelic}.
De fa\c con pr\'ecise, leur \'enonc\'e est le cas particulier de l'\'enonc\'e ci-dessous
 lorsque le morphisme $X \to Y$  provient d'un corps global $K$ dont 
les $K_{v}$ sont les com\-pl\'e\-t\'es.

\medskip

\noindent{\bf Th\'eor\`eme  5.7.} 
{\it Soient $(K_{v})_{v\in V}$ et $f:X\to Y$
 comme dans \rref{ssecPPA}. On suppose que:
\begin{itemize}
\item $f$ est universellement ouvert, surjectif, \`a fibres g\'eom\'etriquement int\`egres;
\item les corps valu\'es $K_{v}$ sont hens\'eliens;
\item pour tout $q\in\NN$ et pour presque tout $v\in V$%
\footnote{Bien entendu, l'ensemble des \og bons\fg\ $v$ d\'epend de $q$.}, le corps $K_{v}$ est 
ul\-tra\-m\'e\-trique, \`a corps r\'esiduel \emph{fini} $k_{v}$ de cardinal $>q$. (En particulier, $V_{\infty}$ est fini).
\end{itemize}
Alors:
\begin{numlist}
\item\label{teo_fibre_integre1} Il existe un sous-ensemble fini $S$ de $V$ et un $S$-mod\`ele 
$F: \gX \to \gY$ de $f$ tel que pour tout $v\in V\setminus S$ l'application induite 
$\gX(O_{v}) \to \gY(O_{v})$ soit surjective. (En particulier, $f$ est presque propre.) 

\item\label{teo_fibre_integre2} Si $f$ est lisse, l'application $f_{\Top}: \gX(\ade_{\soul{K}})_{\Top} \to \gY(\ade_{\soul{K}})_{\Top}$ est  ouverte.
Plus pr\'ecis\'ement, son image $I$ est ouverte dans $\gY(\ade_{\soul{K}})_{\Top}$, et l'application induite $\ol{f}_{\Top}: \gX(\ade_{\soul{K}})_{\Top} \to I$
admet assez de sections locales, et assez de sections si $V_{\infty}=\emptyset$.

\end{numlist}
}

Nous convenons d'appeler \emph{ad\'eliquement surjectif} la propri\'et\'e (1) 
du th\'eo\-r\`eme.
L'in\-gr\'e\-dient principal est l'estim\'ee suivante de type \og Lang-Weil \fg.


\medskip

\noindent{\bf Lemme 5.8.} 
{\it Soit $A$ un anneau et 
soit $F: E \to B$ un $A$-morphisme entre $A$-sch\'emas de pr\'esentation  finie. Soit $d$ un entier $\geq 1$.

\begin{numlist}
\item \label{lem_LW1}
 Il existe une constante
 $C_1>0$ telle que pour tout corps fini $k$ de cardinal $q$ 
 qui est une $A$-alg\`ebre et pour tout point $b \in B(k)$
 tel que $E_b$ soit de dimension $d$, on a 
 $$
\mid  \# E_{b}( k ) \mid \, \leq \, \, \,  C_1 \, q^{ d}.
 $$

\item \label{lem_LW2} Il existe une constante
 $C_2>0$ telle que pour tout corps fini $k$ de cardinal $q$ 
 qui est une $A$-alg\`ebre et pour tout point $b \in B(k)$
 tel que $E_b$ soit g\'eom\'etriquement  {irr\'eductible} de dimension $d$, on a 
 $$
\mid  \# E_{b}( k ) \, - \, q^d  \mid \, \leq \, \, \,  C_2 \, q^{ d - \frac{1}{2}}.
 $$
 
 \item \label{lem_LW3} Il existe une constante $C_3>0$ telle que pour tout corps fini 
 $k$ de cardinal $q \geq  C_3$ 
 qui est une $A$-alg\`ebre,  et pour tout point $b \in B(k)$
 tel que $E_b$ soit 
 {g\'eom\'etriquement int\`egre}, 
 alors $E_{b}$ a un $k$-point lisse.  
\end{numlist}
}

\medskip

\noindent{\bf Remarque 5.9.} 
{\rm 
(a) Pour l'application au th\'eor\`eme 5.7, on
n'a besoin que 
du cas lisse et s\'epar\'e qui est trait\'e sur $\ZZ$ par Conrad
\cite[Lemma 4.6]{ConAdelic}, toujours dans le cas o\`u le point 
$b$ est de corps r\'esiduel de cardinal $q$. 

\smallskip

\noindent (b) Pour le th\'eor\`eme 5.7 et plus loin (corollaire \ 6.7), 
seule la partie (3) de l'\'enonc\'e sera utilis\'ee. 
Un cas particulier important est celui
o\`u $F: E \to B $ est un espace homog\`ene sous un $B$-sch\'ema en groupes
$G$ lisse, s\'epar\'e, de pr\'esentation finie, et \`a fibres 
g\'eom\'etriquement connexes. De ce point de vue,
ce cas particulier est bien plus simple. En effet, 
d'apr\`es Lang \cite[Th. 2]{L}, pour tout corps fini $k$ 
 qui est une $A$-alg\`ebre et pour tout point $b \in B(k)$,
 on a $E_{b}( k ) \not =\emptyset$. 
}

\smallskip

La version ci-dessus est une l\'eg\`ere g\'en\'eralisation
de celle de  Poonen \cite[Th. 7.7.1]{Po}.

\medskip

{\sl \noindent D\'emonstration du lemme 5.8 \hskip0.1em: } si $A=\ZZ$ et si  on se restreint aux points 
$b$ dont le corps r\'esiduel est de cardinal $q$, il s'agit 
exactement du r\'esultat cit\'e. 
Par inspection de la d\'emonstration, celle-ci se g\'en\'eralise au cas d'un 
point arbitraire $b \in B(k)$.

Dans le cas g\'en\'eral, il existe un sous-anneau $A_{0}$ de $A$ qui est une $\ZZ$-alg\`ebre de type fini et un $A_0$-morphisme 
$F_0: E_0 \to B_0$ entre $A_0$-sch\'emas de pr\'esentation  finie tel
que $F= F_0 \times_{A_0} A$. Comme $F_{0}$ est en particulier un morphisme de $\ZZ$-sch\'emas de pr\'esentation finie, le cas $A=\ZZ$ s'applique.
\qed
\medskip

Pour la preuve du th\'eor\`eme 5.7, nous aurons aussi besoin du r\'esultat suivant: 

\medskip

\noindent{\bf Lemme 5.10.} 
{\it Soit $f:X\to Y$ un morphisme de sch\'emas. On suppose que $Y$ est r\'eduit et que $f$ est localement de pr\'esentation finie, universellement ouvert, \`a fibres g\'eo\-m\'e\-tri\-que\-ment r\'eduites. Alors $f$ est plat.
}

\smallskip

\dem lorsque $Y$ est localement noeth\'erien cela r\'esulte directement de \cite[(15.2.3)]{EGA4}. En g\'en\'eral, on peut supposer $X$ et $Y=\Spec(A)$ affines; il existe un sous-anneau $A_{0}$ de $A$, de type fini sur $\ZZ$, tel que $f$ provienne par changement de base d'un morphisme de type fini $f_{0}:X_{0}\to Y_{0}=\Spec(A_{0})$. D'apr\`es  \cite[th. 6.6]{R}, on peut choisir  $A_{0}$ de sorte que $f_{0}$ soit universellement ouvert. D'autre part, l'ensemble $E$ des $y\in Y_{0}$ tel que $f_{0}^{-1}(y)$ soit g\'eom\'etriquement r\'eduit est constructible dans $Y_{0}$ \cite[(9.7.7)]{EGA4} et l'hypoth\`ese sur $f$ assure que l'image de $Y$ dans $Y_{0}$ est contenue dans $E$. On en d\'eduit par \cite[(8.3.4)]{EGA4} qu'en choisissant $A_{0}$ assez grand, on a $E=Y_{0}$. Il r\'esulte alors du cas noeth\'erien que $f_{0}$, et donc $f$, est plat.\qed

\medskip
{\sl \noindent D\'emonstration du th\'eor\`eme 5.7 \hskip0.1em: } \ref{teo_fibre_integre1} En \'ecrivant $\ade_{\soul{K}}$ comme  limite inductive 
des $\ade_{\soul{K},S}$, il existe $S_{\sharp}\subset V$ fini et 
un $S_\sharp$-mod\`ele $F_\sharp: \gX_\sharp \to \gY_\sharp$ de $X \to Y$ 
(avec $\gX_\sharp$ et $\gY_\sharp$ s\'epar\'es de pr\'esentation finie sur
$\ade_{\soul{K},S_\sharp}$). En outre, quitte \`a \'etendre $S_\sharp$,  on peut 
supposer $F_\sharp$ universellement ouvert d'apr\`es \cite[th. 6.6]{R}, surjectif  \cite[(8.10.5)]{EGA4}, et de plus \`a fibres g\'eom\'etriquement int\`egres: pour  ce dernier point on observe que  pour $S$ fini con\-te\-nant $S_{\sharp}$, l'ensemble des points $y\in \gY_{S}$ tels que $F_{S}^{-1}(y)$ soit g\'eom\'etriquement int\`egre  est constructible dans $\gY_{S}$ \cite[9.7.7]{EGA4}, et l'on conclut en utilisant \cite[(8.3.4)]{EGA4} comme dans la preuve de 5.10.

Par application du lemme 5.8.\ref{lem_LW3},  il existe une constante $C>0$ telle que
pour tout corps fini $k$ de cardinal $q>C$ qui est une $\ade_{\soul{K},S_\sharp}$-alg\`ebre et pour tout point $y \in \gY_\sharp(k)$, on a 
$ \gX_{\sharp,y}^{\mathrm{lisse}}(k)  \not= \emptyset$, o\`u $\gX_{\sharp,y}^{\mathrm{lisse}}$ d\'esigne le
lieu lisse du $k$-sch\'ema $\gX_{\sharp,y}$.
On pose alors  
$$S= S_\sharp \cup  V_{\infty} \cup \bigl\{ v \in V \setminus V_\infty \, 
\mid \, \# k_v   \leq C  \bigr\}.$$
C'est un ensemble fini; notons
 $F:\gX\to\gY$ le 
$S$-mod\`ele de $f$ d\'eduit de $F_\sharp$, et v\'erifions que la propri\'et\'e \ref{teo_fibre_integre1} est satisfaite. Soient donc $v\in V\setminus S$ et  $y\in \gY(O_{v})$. Le point $\ol{y}\in\gY(k_{v})$ d\'eduit de $y$ se rel\`eve (puisque $\# (k_{v})>C$)
en un point $\ol{x}\in\gX(k_{v})$, lisse dans sa fibre. Consid\'erons $F_{\red}:\gX_{\red}\to\gY_{\red}$, restriction de $F$ au-dessus de $\gY_{\red}$. Alors  il r\'esulte de 5.10 que $F_{\red}$ est plat, et donc lisse au point $\ol{x}$. Comme $y:\Spec(O_{v})\to \gY$ se factorise par $\gY_{\red}$, et que $O_{v}$ est hens\'elien, on conclut qu'il existe $x:\Spec(O_{v})\to \gX_{\red}$ relevant $y$ et prolongeant $\ol{x}$.

\smallskip

\noindent \ref{teo_fibre_integre2} est cons\'equence de  \ref{teo_fibre_integre1} et des assertions \ref{prop_pp1_1} et  \ref{prop_pp1_3} de 5.5.\,\ref{prop_pp1}.
\qed

\medskip

\noindent{\bf Remarque 5.11.}
Les conditions sur $\soul{K}$ de 5.7 sont v\'erifi\'ees lorsque $K$
est un corps global, $V$ l'ensemble de ses places et $(K_{v})$ la famille de ses compl\'et\'es; c'est le cas consid\'er\'e dans 
\cite{O} et dans \cite{ConAdelic}.

\section{Le cas des torseurs}\label{sec:torseurs}
\subsection{Rappels et conventions}
\subsubsection{Torseurs. }\label{ssec:convtorseurs}Soient $S$ un sch\'ema et $H$ un 
$S$-sch\'ema en groupes. Nous entendrons par \emph{$H$-torseur} sur $S$ un faisceau 
sur $S$ pour la topologie fppf, muni d'une action de $H$ (ou plut\^ot du faisceau qu'il
repr\'esente) et qui est un $H$-torseur pour cette action. Lorsque $H$ est de 
pr\'esentation finie sur $S$,
un tel torseur $X$ est automatiquement un $S$-espace alg\'ebrique  
de pr\'esentation finie \cite[78.11.8, tag 04U1]{St},
s\'epar\'e si $H$ l'est; nous dirons que $X$ est \emph{repr\'esentable} si c'est un sch\'ema, ce qui est toujours le cas si $H$ est affine sur $S$, 
ou encore si $S$ est le spectre d'un corps. 

Pour $S$ et $H$ comme ci-dessus, si $Y$ est un $S$-sch\'ema, 
un $H\times_{S}Y$-torseur (repr\'esentable) sera aussi  appel\'e, par abus, un \og $H$-torseur (repr\'esentable) au-dessus de $Y$\fg.

\subsubsection{Cas des ad\`eles. }Dans la suite du \S\ref{sec:torseurs}, on fixe 
un ensemble d'indices $V$ et 
une famille  $\soul{K}=(K_{v},\abs_{v})_{v\in V}$ de
corps valu\'es. 

Si $G$ est un $\ade_{\soul{K}}$-sch\'ema en groupes s\'epar\'e de pr\'esentation finie,
$Y$ un $\ade_{\soul{K}}$-sch\'ema s\'epar\'e de pr\'esentation finie, 
et $f:X\to Y$ un $G$-torseur repr\'esentable au-dessus de $Y$, on se propose d'\'etudier 
l'application continue
$$f_{\Top}:X(\ade_{\soul{K}})_{\Top}\to Y(\ade_{\soul{K}})_{\Top}$$
induite par $f$. Il est imm\'ediat que le groupe topologique $G(\ade_{\soul{K}})_{\Top}$ op\`ere librement sur $X(\ade_{\soul{K}})_{\Top}$, et que les fibres non vides de $f_{\Top}$ sont les orbites pour cette action. De fa\c{c}on plus pr\'ecise, $f_{\Top}$ est un \emph{pseudo-torseur topologique} sous $G(\ade_{\soul{K}})_{\Top}$, c'est-\`a-dire que l'application 
$$\begin{array}{ccc}
G(\ade_{\soul{K}})_{\Top}\times_{Y(\ade_{\soul{K}})_{\Top}}X(\ade_{\soul{K}})_{\Top} & \ffl & X(\ade_{\soul{K}})_{\Top}\times_{Y(\ade_{\soul{K}})_{\Top}}X(\ade_{\soul{K}})_{\Top}\\
(g,x) & \longmapsto & (gx,x)
\end{array}$$
est un hom\'eomorphisme; c'est l\`a une cons\'equence formelle 
de la compatibilit\'e aux produits fibr\'es.

\subsection{Torseurs presque propres: une condition suffisante}

Soit $G$ un $\ade_{\soul{K}}$-sch\'ema en groupes s\'epar\'e de pr\'esentation finie. Lorsque nous parlerons d'un $S$-mod\`ele $\gG$ de $G$, pour une partie finie $S$ de $V$, nous sup\-po\-se\-rons toujours $\gG$ muni d'une structure de $\ade_{\soul{K},S}$-sch\'ema 
en groupes in\-dui\-sant celle de $G$. Nous nous int\'eressons ici \`a la condition suivante sur $G$:
\smallskip

\noindent NT($G$):  \hfill \begin{minipage}[t]{12cm} il existe une partie finie 
$S$ de $V$ et un $S$-mod\`ele $\gG$ de $G$ tels que, pour presque tout $v\in V\setminus S$,
l'application naturelle
$$\rho(\gG,v):\qquad {\rH}^1(O_v,\gG) \to {\rH}^1(K_v,\gG)$$
ait un noyau trivial.
\end{minipage}
\medskip

\noindent Noter que la condition en question 
est alors satisfaite par \emph{tout} mod\`ele
de $G$.

\medskip

\noindent{\bf Proposition 6.1.} 
{\it  Soit $G$ un $\ade_{\soul{K}}$-sch\'ema en groupes 
s\'epar\'e de pr\'e\-sen\-ta\-tion finie, $Y$ un $\ade_{\soul{K}}$-sch\'ema s\'epar\'e
de pr\'esentation finie et $f: X \to Y$ un $G$-torseur repr\'esentable.  On suppose que  
\textup{NT}($G$) est v\'erifi\'ee. Alors $f$ est presque propre.
}

\smallskip

\dem  on peut choisir un $S$-mod\`ele $\gG$ de $G$ tel que  l'application 
$\rho(\gG,v)$ ait un noyau trivial pour \emph{tout} $v\notin S$, et que 
$f$ admette un $S$-mod\`ele $F:\gX\to \gY$ qui soit un $\gG$-torseur
en vertu de \cite[VI$_B$.10.16]{SGA3}. 

Fixons $v\notin S$ et notons $F_{v,\Top}:\gX(O_{v})\to\gY(O_{v})$ l'application 
induite par $F$. Soit $y\in \gY(O_{v})$. On a alors les \'equivalences:
$$\begin{array}{ccl}
y\in \im(F_{v,\Top})&
\Leftrightarrow& \text{le $\gG_{v}$-torseur $y^*\gX$ sur $\Spec\;(O_{v})$ est trivial}\\
&\Leftrightarrow& \text{le $\gG_{K_{v}}$-torseur $y_{K_{v}}^*\gX$ sur $\Spec\;(K_{v})$ est trivial}\\
&\Leftrightarrow& y_{K_{v}}\in \im\left(\gX(K_{v})\to\gY(K_{v})\right)
\end{array}$$
o\`u la seconde \'equivalence r\'esulte de l'hypoth\`ese sur  $\rho(\gG,v)$. Le
lemme en r\'esulte.\qed

\medskip

\noindent{\bf Remarque 6.2.} {\rm La condition NT($G$) est satisfaite si $G$ est propre sur $S$.
On retrouve ainsi que $f$ est presque propre dans ce cas (remarque 5.3\,\ref{rem_propre2}).
}

\medskip

La proposition 5.5 a la cons\'equence suivante.

\medskip

\noindent{\bf Corollaire 6.3.} 
{\it Sous les hypoth\`eses de 6.1 (incluant $\textup{NT}(G)$), on suppose de 
plus que $G$ est \emph{lisse} et que
 tous les corps valu\'es $K_v$ sont hens\'eliens.
On note $I$ l'image de $f_\Top: X(\ade_{\soul{K}})_\Top \to Y(\ade_{\soul{K}})_\Top$. 
Alors:

\begin{numlist}
\item   $I$ est ferm\'e dans $Y(\ade_{\soul{K}})_\Top$.
\item  Si de plus $V_\infty$ est fini, l'application induite 
$X(\ade_{\soul{K}})_\Top \to I$ est une 
$G(\ade_{\soul{K}})_\Top$-fibration principale, triviale si $V_{\infty}=\emptyset$.
\end{numlist}
}

\smallskip

\dem La proposition 6.1 montre que $f$ est presque propre.
\smallskip

\noindent (1) Comme $G$ est lisse, $f$ est lisse et les applications
$X(K_v) \to Y(K_v)$ sont ouvertes (prop.\ 3.5) et en 
particulier d'image $I_v$ ouverte.
Pour chaque $v$, le compl\'ementaire de $I_v$
est une r\'eunion d'images de torseurs sur $Y_{K_v}$ sous des formes tordues de $G$, 
donc ce compl\'ementaire est ouvert (l'argument est d\'etaill\'e dans \cite[prop. 3.4.1]{GGMB}). Ainsi $I_v$ est ouvert et ferm\'e pour tout $v\in V$.
La proposition 5.5\,\ref{prop_pp2} montre que l'image de $f_\Top$ est ferm\'ee.
\smallskip

\noindent (2) La   proposition 5.5\,\ref{prop_pp1}\,\ref{prop_pp1_1} montre que l'application induite 
$X_\Top \to I$ a assez de sections locales, c'est donc une 
$G(\ade_{\soul{K}})$-fibration principale. Si $V_{\infty}=\emptyset$, elle a m\^eme assez de sections, donc est triviale.
\qed
\medskip

\subsection{Quelques cas o\`u la propri\'et\'e d'injectivit\'e est connue}

Soient $A$ un anneau int\`egre, $F$ son corps de fractions.
Soit $H$ un $A$-sch\'ema en groupes plat de pr\'e\-sen\-ta\-tion finie.
On sait que l'application ${\rH}^1(A,H) \to {\rH}^1(F,H)$ est injective dans les cas suivants:

\begin{numlist}

\item $H$ est propre sur $A$ et $A$ est un anneau de Pr\"ufer (ses anneaux locaux sont
des anneaux de valuation);
\item $H$ est fini sur $A$ et $A$ est normal (dans ce cas et dans le pr\'ec\'edent, on 
a m\^eme $X(A)=X(F)$ pour tout $H$-torseur $X$ sur $\Spec(A)$);
\item $A$ est semi-local r\'egulier et $H$ est de type multiplicatif \cite[th. 4.1]{CTS};
 \item\label{cas-injectifs4} 
$A$ est un anneau de valuation et $H$ est
 r\'eductif (Bruhat-Tits pour le cas d'un anneau de valuation discr\`ete complet, voir \cite[Th. I.1.2.2]{G1}; 
 voir \cite{N} et \cite[Theorem 1.3]{Guo} pour le cas g\'en\'eral).

 \item $A=k[[t]]$, o\`u $k$ est un corps, $H$ est lisse et provient de $k$ \cite[th. 5.5]{FG}.
\end{numlist}

On peut attraper d'autres cas via le d\'evissage facile suivant.

\medskip

\noindent{\bf Lemme 6.4.} 
{\it Soient $A$ un anneau int\`egre, $F$
son corps de fractions.
Soit $1 \to G_1 \to G_2 \to G_3 \to 1$ une suite exacte de 
$A$-sch\'emas en groupes. 
On suppose que les deux applications ${\rH}^1(A,G_1) \to {\rH}^1(F,G_1)$ et 
${\rH}^1(A,G_3) \to {\rH}^1(F,G_3)$ ont un noyau trivial. Alors on a aussi
$$  \ker\bigl({\rH}^1(A,G_2) \to {\rH}^1(F,G_2) \bigr)=1$$
dans les deux cas suivants:
\begin{numlist} 
 \item\label{lem_propre1} l'application naturelle  $G_3(A)\to G_3(F)$ est surjective;
 \item\label{lem_propre2}  ${\rH}^1(A,G_1)=1$.
 \end{numlist}
}

\smallskip

\dem la suite de l'\'enonc\'e donne lieu au diagramme commutatif exact 
d'ensembles point\'es \cite[III.3.3]{Gi}
$$
\begin{CD}
 G_3(A) @>{\varphi_A}>> {\rH}^1(A,G_1)  @>{\alpha_A}>> {\rH}^1(A,G_2) @>{\beta_A}>> {\rH}^1(A,G_3) \\
@VV{}V  @V{a_{G_1}}VV @V{a_{G_2}}VV @V{a_{G_3}}VV\\
 G_3(F) @>{\varphi_F}>> {\rH}^1(F,G_1) 
 @>{\alpha_F}>> {\rH}^1(F,G_2) @>{\beta_K}>> {\rH}^1(F,G_3) .\\
\end{CD}
$$
Le cas \ref{lem_propre2} est alors imm\'ediat et laiss\'e au lecteur. Dans le cas 
\ref{lem_propre1}, on part d'une classe $\gamma_2 \in \ker a_{G_{2}}$. 
Comme $\ker a_{G_{3}}=1$ on voit qu'il existe  $\gamma_{1}\in {\rH}^1(A,G_1)$ s'appliquant
sur $\gamma_{2}$ et telle que $\gamma_{1,F}:=a_{G_{1}}(\gamma_{1})$ soit dans le 
noyau de $\alpha_{F}$. 

Le groupe $G_3(A)$ agit  \`a droite sur ${\rH}^1(A,G_1)$ et  
les fibres de $\alpha_A$ co\" {\i}n\-cident avec les orbites pour cette action \cite[III.3.3.3]{Gi}.
Il en est de m\^eme pour $G_3(F)$ et $\alpha_F$.
Ainsi il existe $g_{3,F} \in G_3(F)$ tel que $\gamma_{1,F} \, . \, g_{3,F} =1  \in {\rH}^1(F,G_1) $. Vu l'hypoth\`ese \ref{lem_propre1}, $g_{3,F}$ se rel\`eve en  $g_{3} \in G_3(A)$ et
l'on a dans ${\rH}^1(F,G_1) $ l'\'egalit\'e
$$1=\gamma_{1,F} \, . \, g_{3,F}=(\gamma_{1} \, . \, g_{3})_{F}$$
d'o\`u $\gamma_{1} \, . \, g_{3}\in \ker a_{G_{1}}=1$. En conclusion, 
$\gamma_{2}=\alpha_{A}(\gamma_{1})=\alpha_{A}(\gamma_{1} \, . \, g_{3})=1$.
\qed
\medskip

Par applications successives de 6.4 et des crit\`eres \'enonc\'es 
juste avant, on en d\'eduit par exemple:

\medskip

\noindent{\bf Corollaire 6.5.} 
{\it On suppose que $A$ est un anneau de valuation hens\'elien, de corps r\'esiduel $\kappa$. Soit  $G$ un $A$-sch\'ema en groupes de pr\'esentation finie  admettant une suite de composition
$$ G_{1} \vartriangleleft G_{2} \vartriangleleft G$$
dans laquelle:
\begin{romlist}
\item $G_{1}$ est lisse et ${\rH}^1(\kappa,G_{1,\kappa})=1$;
\item $G_{2}/G_{1}$ est r\'eductif;
\item $G/G_{2}$ est propre sur $\Spec(A)$.
\end{romlist}
Alors l'application naturelle ${\rH}^1(A,G)\to{\rH}^1(F,G)$ est \`a noyau trivial.
}

\smallskip

\dem (i) entra\^{\i}ne que $\rH^1(A,G_1)=1$.
On consid\`ere la suite exacte $1 \to G_1 \to G_2 \to G_2/G_1 \to 1$.
Si $G_1=G_2$, on a $\rH^1(A,G_2)=1$.
Si $G_2/G_1$ est r\'eductif, alors $\rH^1(A,G_2/G_1) \to \rH^1(F,G_2/G_1)$ est \`a noyau trivial
d'apr\`es le cas \ref{cas-injectifs4} ci-dessus; il suit par d\'evissage que 
$\rH^1(A,G_2) \to \rH^1(F,G_2)$ est \`a noyau trivial.
Ainsi dans tous les cas $\rH^1(A,G_2) \to \rH^1(F,G_2)$ est \`a noyau trivial. 
Comme $G/G_2$ est suppos\'e propre, on a  $(G/G_2)(A)= (G/G_2)(F)$ 
en vertu du crit\`ere valuatif de propret\'e.
Le lemme 6.4.\,\ref{lem_propre2}
permet de conclure que  $\rH^1(A,G) \to \rH^1(F,G)$ est  \`a noyau trivial. 
\qed

\subsection{Applications aux torseurs ad\'eliques}\label{AppTorsAd}

On suppose dans cette section que tous les $K_{v}$ sont \emph{hens\'eliens}, et 
que $V_{\infty}$ est \emph{fini}. 

\medskip

\noindent{\bf Corollaire 6.6.} 
{\it  Soit $G$ un $\ade_{\soul{K}}$-sch\'ema en groupes s\'epar\'e lisse de pr\'esentation finie.
Soient  $Y$ un $\ade_{\soul{K}}$-sch\'ema s\'epar\'e  de pr\'esentation finie
et $f: X \to Y$ un $G$-torseur repr\'esentable. On note $I$ l'image de $f_\Top: X(\ade_{\soul{K}})_\Top \to Y(\ade_{\soul{K}})_\Top$.

On suppose qu'il existe un $S$-mod\`ele
lisse $\gG$ de $G$ qui admet une suite de composition
$$ \gG_{1} \vartriangleleft \gG_{2} \vartriangleleft \gG$$
dans laquelle $\gG_1$ est affine lisse \`a fibres unipotentes d\'eploy\'ees,
$\gG_2/\gG_1$ est r\'eductif et $\gG/\gG_2$ est propre.

\begin{numlist}

\item   $I$ est ferm\'e dans $Y(\ade_{\soul{K}})_\Top$.

\item  L'application induite $X(\ade_{\soul{K}})_\Top \to I$ est une 
 $G(\ade_{\soul{K}})_\Top$-fibration principale.
\end{numlist}
 }

\smallskip

\dem Le corollaire 6.5 garantit que la 
propri\'et\'e $\textup{NT}(G)$ vaut. Ainsi ce corollaire est une cons\'equence
du corollaire 6.3.
\qed

\medskip

Rappelons qu'un corps $k$ est dit \emph{pseudo-alg\'ebriquement clos} si toute $k$-vari\'et\'e g\'eo\-m\'e\-tri\-que\-ment int\`egre 
admet un point rationnel. C'est le cas, entre autres, des corps s\'e\-pa\-ra\-ble\-ment clos et des extensions alg\'ebriques infinies d'un corps fini.

\newpage

\noindent{\bf Corollaire 6.7.} 
{\it Soit $\soul{K}=(K_{v})_{v\in V}$ comme en \rref{ssecProdRestr}.
On suppose que:
\begin{itemize}
\item les corps valu\'es $K_{v}$ sont hens\'eliens;
\item pour chaque  $v\in V$  ultram\'etrique, le corps r\'esiduel
$k_v$ est fini ou pseudo-alg\'e\-bri\-que\-ment clos.
\end{itemize}

\noindent Soit $G$ un $\ade_{\soul{K}}$-sch\'ema en groupes
s\'epar\'e lisse de pr\'esentation finie, \`a fibres connexes.
Soient  $Y$ un $\ade_{\soul{K}}$-sch\'ema s\'epar\'e  de pr\'esentation finie
et $f: X \to Y$ un $G$-torseur. Soit $I\subset Y(\ade_{\soul{K}})$ l'image de $f_\Top$. 
Alors:
\begin{numlist}
\item   $I$ est ferm\'e dans $Y(\ade_{\soul{K}})_\Top$.
\item  Si de plus $V_\infty$ est fini, l'application induite 
$X(\ade_{\soul{K}})_\Top \to I$ est une 
$G(\ade_{\soul{K}})_\Top$-fibration principale, triviale si $V_{\infty}=\emptyset$.
\end{numlist}
}

\smallskip

\dem Il suffit de voir que la condition $\textup{NT}(G)$ est v\'erifi\'ee.
On peut choisir un $S$-mod\`ele $\gG$ de pr\'esentation finie de $G$  tel que  
$f$ admette un $S$-mod\`ele $F:\gX\to \gY$ qui soit un $\gG$-torseur
en vertu de \cite[VI$_B$.10.3 et 10.16]{SGA3}. 
Quitte \`a agrandir $S$, on peut supposer que $\gG$ est s\'epar\'e et lisse
\cite[8.10.5.(ii) et 17.7.8.(2)]{EGA4} et l'argument de la preuve
du th\'eor\`eme 5.7 permet de supposer
$\gG$ \`a fibres connexes. 

Soit $v\in V\setminus(S\cup V_{\infty})$.
On a une bijection $\rH^1(O_v,\gG) \to \rH^1(k_v,\gG)$ en vertu
de \cite[XXIV.8.1]{SGA3}. Nous allons montrer que 
$\rH^1(k_v,\gG)=1$ et partant que $\rH^1(O_v,\gG)=1$. 
Puisque $\gG_{k_v}$ est lisse et connexe, il s'agit du th\'eor\`eme de Lang si $k_v$ 
est  un corps fini \cite[Th. 2]{L}; d'autre part cela entra\^{\i}ne que tout $\gG_{k_v}$-torseur  
est g\'eom\'etriquement int\`egre, donc l'assertion est triviale si $k_{v}$ est pseudo-alg\'e\-bri\-que\-ment clos. 

Le crit\`ere $\textup{NT}(G)$ est donc v\'erifi\'e et le corollaire 6.3 s'applique.
\qed

\subsection{Le cas de caract\'eristique nulle 
}\label{ssec-carac-nulle2}

On fixe un corps $K$ de caract\'eristique nulle  et une famille non vide $\soul{K}=(K_{v},\abs_{v})_{v\in V}$ d'extensions valu\'ees de $K$.
On suppose en outre que:
\begin{itemize}
\item $V_{\infty}$ est fini; 
\item pour tout $v\in V$, $K_{v}$ est hens\'elien;
\item pour tout $f\in K^\times$, on a $\vert f\vert_{v}=1$ pour presque tout $v\in V$.
\end{itemize}
\smallskip

\noindent En particulier,  $K\subset \ade_{\soul{K}}$ d'apr\`es la troisi\`eme condition.

Un cas bien connu est celui d'un corps de nombres muni de la famille de ses compl\'et\'es; un autre est celui o\`u $K$ est le corps des fonctions rationnelles d'un $\QQ$-sch\'ema int\`egre, noeth\'erien et normal  $Z$, $V$ la famille des valuations divisorielles de $Z$, et $K_{v}$ le compl\'et\'e (ou le hens\'elis\'e) de $K$ en $v$.

\newpage

\noindent{\bf Th\'eor\`eme 6.8.} 
{\it  Sous les hypoth\`eses pr\'ec\'edentes, soient $G$ un $K$-sch\'ema en groupes de type fini, $Y$ un $\ade_{\soul{K}}$-sch\'ema s\'epar\'e de pr\'esentation finie et $f:X\to Y$ un $G$-torseur. 
On note $I$ l'image de $f_\Top: X(\ade_{\soul{K}})_\Top \to Y(\ade_{\soul{K}})_\Top$. Alors:
\begin{numlist}
\item   $I$ est ferm\'e dans $Y(\ade_{\soul{K}})_\Top$.
\item  L'application induite $X(\ade_{\soul{K}})_\Top \to I$ est une $G(\ade_{\soul{K}})_\Top$-fibration principale.
\end{numlist}
}

\smallskip

\dem Puisque $K$ est de caract\'eristique nulle, il existe une suite de composition $G_{1} \vartriangleleft G_{2} \vartriangleleft G$ o\`u $G_{1}$ est unipotent d\'eploy\'e, $G_{2}/G_{1}$ est r\'eductif, et $G/G_{2}$ est propre sur $K$. D'apr\`es \cite[\S\ 8]{EGA4}, il existe un sous-anneau $A_{0}$ de $K$, de type fini sur $\ZZ$, tel que $G$ provienne par changement de base d'un $A_{0}$-sch\'ema en groupes lisse admettant une suite de composition (\`a composantes lisses) avec les m\^emes propri\'et\'es. Comme $A_{0}\subset \ade_{\soul{K},S}$ pour $S\subset V$ fini assez grand, on en d\'eduit un $S$-mod\`ele $\gG$ de $G$ v\'erifiant la condition du corollaire 6.6; ce dernier donne le r\'esultat.\qed

\subsection{Exemple d'une image non ferm\'ee et d'une 
application non stricte}
Soit $K$ le corps global $\FF_p(t)$, muni de la famille de ses compl\'et\'es $(K_{v})_{v\in V}$, o\`u $V$ s'identifie \`a l'ensemble des points ferm\'es de $\PP^1_{\FF_{p}}$. On va donner un exemple de $H$-torseur $f: X \to Y$  tel que $f_\Top$ n'est pas stricte et telle que l'image $I$ de $f_\Top$ n'est pas ferm\'ee.
Cet exemple est un avatar global d'un exemple local \cite[\S 7.1]{GGMB}. 

On note $Y_0 =\Aa^1_K$. Le groupe $G_0= \GG_a \rtimes \GG_m$ agit sur $Y_0$ par $(x,y).z= x^p+y^p z$.
On note $H_0$ le stabilisateur de $t \in \Aa^1_K(K)$ pour cette action. Le morphisme d'orbite en $t$ donne un isomorphisme
de $G_0$-vari\'et\'es  $G_0/H_{0} \simlgr \Aa^1_K= Y_0$. 
On pose  $X={G_0 \times_K \ade_{\soul{K}}}$, $Y=Y_0 \times_K \ade_{\soul{K}}$,
$H=H_0 \times_K \ade_{\soul{K}}$,
On note  $f: X \to G/H \simlgr Y$ le morphisme quotient. 

Notant $\ide_{\soul{K}}=\GG_{m}(\ade_{\soul{K}})_{\Top}$ le groupe topologique des id\`eles de $K$,
l'image $I$ de $f_\Top: \ade_{\soul{K}} \rtimes \ide_{\soul{K}} \to \ade_{\soul{K}}$ est  $\bigl\{ x^p+y^p t  \, \mid x \in \ade_{\soul{K}}, \, y \in \ide_{\soul{K}} \bigr\}$.
Alors $0$ n'appartient pas \`a $I$ mais est adh\'erent \`a $I$ donc $I$ n'est pas ferm\'ee.

Observons que $f_\Top$ est injective (de fa\c{c}on \'equivalente, $H_{\Top}$ est trivial): en effet, pour tout $v\in V$, $t$ n'est pas une puissance $p$-i\`eme dans $K_{v}$. Supposons,  par l'absurde, que $f_\Top:\ade_{\soul{K}} \rtimes \ide_{\soul{K}} \to \ade_{\soul{K}}$ soit  stricte.
Comme elle est injective, c'est donc un plongement topologique.
En particulier, le compos\'e $\ide_{\soul{K}} \hookrightarrow {\ade_{\soul{K}} \rtimes \ide_{\soul{K}}} \, \xrightarrow{ \,   f_\Top \, } \,  \ade_{\soul{K}}$,
$y \to y^p t$  est un plongement  topologique, ce qui est absurde. Ainsi $f_\Top$ n'est pas stricte et 
n'induit pas une $H_\Top$-fibration topologique. 

Enfin, $I$ est localement ferm\'e dans $\ade_{\soul{K}}$:  il suffit de le v\'erifier au voisinage de $t$.
 On consid\`ere l'ouvert $\Omega= t+\ade_{\soul{K}, \emptyset}
= t + \prod_{v \in V} O_v$ de $\ade_{\soul{K}}$.
Vu que $t+ O_v \cap 
f_v(K_v \rtimes K_v^\times)$ est ferm\'e dans $t+O_v$ pour tout 
$v \in V$ le lemme 2.2.\ref{lemTopProduit3}
indique que $I \cap \Omega$ est ferm\'e dans $\Omega$.

\section{Appendice: ultraparacompacit\'e}

\bigskip

On rappelle qu'un espace topologique $X$
est \emph{ultraparacompact} s'il est s\'epar\'e et 
si  tout recouvrement ouvert $(U_i)_{i \in I}$ de $X$
admet un raffinement ouvert $(V_j)_{j \in J}$ tel que
$X= \coprod_{j \in J} V_j$.

\medskip

\noindent{\bf Proposition 7.1.} 
{\it Soit $(X_{i})_{i\in\FI}$ un syst\`eme projectif d'ensembles, dont l'ensemble d'indices $\FI$ est \emph{totalement ordonn\'e}. Notons $X$ l'espace topologique $\varprojlim_{i\in\FI}X_{i}$, chaque $X_{i}$ \'etant muni de la topologie discr\`ete. Alors tout sous-espace de $X$ est ultraparacompact.
}

\smallskip

\dem quitte \`a remplacer  $\FI$ par un sous-ensemble cofinal, nous pouvons le supposer bien ordonn\'e. Pour $i\in\FI$, notons $X\xrightarrow{p_{i}}X_{i}$ 
l'application naturelle. L'ensemble des $p_{i}^{-1}(b)$, pour $i\in \FI$ et $b\in X_{i}$, est une base d'ouverts de $X$, qu'il est commode de voir comme des \og boules\fg\ dont l'ensemble des rayons est $\FI^\circ$ (l'ensemble totalement ordonn\'e oppos\'e \`a $\FI$): pour $x\in X$, nous poserons donc $\cB(x,i):=p_{i}^{-1}(p_{i}(x))$. Pour $i\leq j$ dans $\FI$ et $x$, $y$ dans $X$, on a les \'equivalences:
$$\begin{array}{lclcl}
\cB(x,i)\cap\cB(y,j)\neq\emptyset &\Leftrightarrow & y\in\cB(x,i)  \\
&\Leftrightarrow & \cB(y,j) \subset\cB(x,i)  & \Leftrightarrow & \cB(y,i) = \cB(x,i).
\end{array}$$

Pour \'etablir la proposition, il est facile de voir \cite[prop. 7]{VN} qu'il suffit de montrer que tout ouvert de $X$ est ultraparacompact. Soit donc $(U_{\lambda})_{\lambda\in L}$ une famille d'ouverts de $X$, de r\'eunion $U$. Nous allons construire une partition de $U$ par des boules dont chacune est contenue dans un $U_{\lambda}$, ce qui  implique le r\'esultat.

Pour chaque $x\in U$, il existe un plus petit indice $i\in\FI$ tel que $\cB(x,i)$ soit contenu dans l'un des $U_{\lambda}$: notons-le $r(x)$ et posons $V(x):= \cB(x,r(x))$. Ainsi $V(x)$ est la plus grande boule contenant $x$ et contenue dans un $U_{\lambda}$. Il est clair que l'ensemble des $V(x)$ recouvre $U$ et raffine $(U_{\lambda})$; d'autre part, si $V(x)\cap V(y)\neq\emptyset$, les \'equivalences \'enonc\'ees plus haut (et la d\'efinition de $r(x)$) impliquent imm\'ediatement que $V(x) = V(y)$, de sorte que l'on a bien une partition de $U$.\qed

\newpage

\noindent{\bf Corollaire 7.2.} 
{\it Soit $A$ un anneau local topologique s\'epar\'e admettant une base $\FI$ de voisinages de $0$ form\'ee d'id\'eaux, et \emph{totalement ordonn\'ee} par inclusion. Soit $\FX$ un $A$-sch\'ema localement de type fini. Alors tout sous-espace de $\FX(A)_{\Top}$  est ultraparacompact.
}

\smallskip

\dem cela r\'esulte de 7.1 et de la description donn\'ee en \ref{ssecLinTopLoc} de $\FX(A)_{\Top}$  comme sous-espace de $\varprojlim_{J\in\FI}\FX(A/J)$.\qed

\medskip

\noindent{\bf Remarque 7.3.}
La condition de l'\'enonc\'e est satisfaite notamment si $A$ est un anneau de valuation (muni de la topologie associ\'ee), ou encore si $A$ est un anneau local muni de la topologie $I$-adique, o\`u $I\subset A$ est un id\'eal tel que $\bigcap_{n\in\NN}I^n=\{0\}$.

\medskip

\noindent{\bf Proposition 7.4.} 
{\it Soit $(K,v)$ un corps valu\'e  ultram\'etrique.  
Soit $X$ un $K$-sch\'ema s\'epar\'e de type fini.
Alors tout sous-espace de $X(K)_{\Top}$ est ultraparacompact.
}

\smallskip

\dem On note $A$ l'anneau de $v$ et $\Gamma$ son groupe. Tout d'abord, l'espace topologique $X(K)_{\Top}$ est s\'epar\'e.
Suivant \cite[(8.8.2)(ii) et (8.10.5)(v)]{EGA4}, il existe $f \in A$ non nul et
un $A_f$-sch\'ema s\'epar\'e de type fini $\gX$ tel que $\gX \times_{A_f} K=X$.
L'application du th\'eor\`eme de compactification de Nagata \cite[thm. 4.1]{ConNagata}
au morphisme $\gX \to \Spec(A_f) \to \Spec(A)$ montre qu'il existe
une $A$-immersion ouverte $\phi: \gX \to \gX^c$ 
o\`u  $\gX^c$ est un $A$-sch\'ema propre. En particulier on a une immersion ouverte $\phi_K: X \to \gX^c_K$. 
 Alors $X(K)$ s'identifie \`a un ouvert de $\gX^c_{K}(K)=\gX^c(K)$ et ce dernier est 
hom\'eomorphe \`a $\gX^c (A)$ en vertu de 
3.4\,\ref{LemChgtAnn2} et de la propret\'e de $\gX^c$; on peut donc conclure par le corollaire 7.2.
\qed
\medskip

\noindent{\bf Remarque 7.5.}
Pour la preuve de 7.4, on peut se passer du th\'eor\`eme de Nagata si l'on sait {\sl a priori} qu'il existe une immersion de $X$ dans la fibre g\'en\'erique d'un $A$-sch\'ema propre; c'est le cas notamment si $X$ est quasi-projectif sur $K$.

\bigskip

\noindent P. {\sc Gille}: \\
C.N.R.S., Institut Camille Jordan - Universit\'e Claude Bernard Lyon 1,
43 boulevard du 11 novembre 1918, 69622 Villeurbanne Cedex (France). \medskip

\noindent L. {\sc Moret-Bailly}: \\
IRMAR, Universit\'{e} de Rennes, Campus de Beaulieu, 35042 Rennes Cedex (France);\\
Centre Henri Lebesgue, programme ANR-11-LABX-0020-0. \medskip

\noindent  Les deux auteurs b\'en\'eficient du soutien du projet Geolie, ANR-15-CE
40-0012 (Agence Nationale de la Recherche, France).


\bigskip

\begin{thebibliography}{EGA IV}

\bibitem{Bt} {\sc B. Bhatt}, {\sl Algebraization and Tannaka duality}, 
 Camb. J. Math. {\bf 4} (2016), 403--461.
%
\bibitem{BTG} {\sc N. Bourbaki}, "Topologie g\'en\'erale, chapitres 1 \`a 4", Springer, 2007.
%
\bibitem{CTS} {\sc J.-L. Colliot-Th\'el\`ene} et {\sc J.-J. Sansuc},  {\sl Principal homogeneous spaces under flasque tori: Applications},
  J. Algebra {\bf 106} (1987), 148--205.
%
\bibitem{ConAdelic} {\sc B. Conrad}, {\sl Weil and Grothendieck Approaches to Adelic Points},  L'Ens. Math. (2) {\bf58} (2012), 61--97.
%
\bibitem{ConNagata} {\sc B. Conrad}, {\sl Deligne's notes on
 Nagata compactifications},  J. Ramanujan Math. Soc. {\bf22} (2007), 205--257.
%
\bibitem{DG} {\sc M. Demazure} et {\sc P. Gabriel},
"Groupes alg\'ebriques", Masson, 1970.
%

\bibitem{SGA3} {\sc M. Demazure} et {\sc A. Grothendieck}, "S\'eminaire de G\'eom\'etrie alg\'ebrique de
l'I.H.E.S., 1963-1964, Sch\'e\-mas en groupes",  
Documents ma\-th\'e\-ma\-tiques vol. 7  et 8, Soci\'et\'e math\'ematique de France, 2011.
%
\bibitem{D} {\sc J. Dugundji}, "Topology", Allyn and Bacon, 1966.
%
\bibitem{E} {\sc O. Endler},  "Valuation Theory", Universitext, Springer, 1972.
%
\bibitem{EP} {\sc A.J. Engler} et {\sc A. Prestel}, "Valued Fields", Springer Monographs in Mathematics, Springer, 2005.
%
\bibitem{FG}  {\sc M. Florence} et  {\sc P. Gille},   {\sl 
Residues on affine Grassmannians},  J. Reine Angew. Math. {\bf 776} (2021), 119--150.
%
\bibitem{GGMB}  {\sc O. Gabber, P. Gille} et {\sc L. Moret-Bailly}, {\sl Fibr\'es principaux sur les corps valu\'es hens\'eliens},
 Algebr. Geom. {\bf 1} (2014), 573--612.
%
\bibitem{G1} {\sc P. Gille}, {\sl R-\'equivalence et torseurs sur la droite affine}, th\`ese de doctorat, Orsay (1994), 
\href{http://math.univ-lyon1.fr/homes-www/gille/divers.html}{http://math.univ-lyon1.fr/homes-www/gille/divers.html}.
%
\bibitem{Gi} {\sc J. Giraud}, "Cohomologie non ab\'elienne",
Grundlehren der mathematischen Wissenschaften {\bf 179},
Springer, 1971.
%
\bibitem{Gre} {\sc M.J. Greenberg}, {\sl Rational points in Henselian discrete valuation rings},  Publ. Math. Inst. Hautes \'Etudes Sci.  {\bf31} (1966), 59--64.
%
\bibitem{EGA1} {\sc A. Grothendieck} et {\sc J.A. Dieudonn\'e}, "\'El\'ements de g\'eom\'etrie alg\'ebrique  I", Grundlehren der mathematischen Wissenschaften vol. 166,
Springer, 1971.
%
\bibitem{EGA4} {\sc A. Grothendieck} et {\sc J.A. Dieudonn\'e}, {\sl
\'El\'ements de g\'eom\'etrie alg\'ebrique 
{\uppercase\expandafter{\romannumeral 4}}.2, 
{\uppercase\expandafter{\romannumeral 4}}.3, et {\uppercase\expandafter{\romannumeral 4}}.4},   Publ. Math. Inst. Hautes \'Etudes Sci. {\bf 24},  {\bf 28} (1965) et {\bf 32} (1967).
%
%
\bibitem{Guo} {\sc N. Guo}, {\sl The Grothendieck-Serre conjecture over valuation rings}, preprint (2020), <https://arxiv.org/abs/2008.02767> (31/05/2022).
%
\bibitem{L} {\sc S. Lang}, {\sl Algebraic groups over finite fields},
 Amer. J. Math. {\bf 78} (1956), 555--563. 
%
\bibitem{Ma} {\sc B. Margaux}, {\sl Passage to the limit in non-abelian \v{C}ech cohomology},
 J. Lie Theory {\bf 17} (2007), 591--596. 
%
\bibitem{MB} {\sc L. Moret-Bailly},  {\sl An extension of Greenberg's theorem to 
general valuation rings},  Manuscripta  Math. {\bf139} (2012) \no1, 153--166.
 %
 \bibitem{N} {\sc Y.A. Nisnevich},  {\sl Espaces homog\`enes principaux rationnellement triviaux 
 et arithm\'etique des sch\'emas en groupes r\'eductifs sur les anneaux de Dedekind},   C. R. Math. Acad. Sci. Paris {\bf 299} (1984), 5--8.
%
\bibitem{O} {\sc J. Oesterl\'e}, {\sl Nombres de Tamagawa et groupes unipotents en caract\'eristique $p>0$}, 
 Invent. Math.  {\bf 78} (1984), 13--88.
 %
\bibitem{Po} {\sc B. Poonen}, "Rational points on varieties", Graduate Studies 
  in Mathematics vol. 186, American Mathematical Society,   2017.
%
 \bibitem{R} {\sc D. Rydh}, {\sl  Submersions and effective descent of \'etale morphisms},
 Bull. Soc. Math. France {\bf 138} (2010), 181-230. 

%
\bibitem{St} {\sc The Stacks project authors}, {\sl The Stacks Project}, \url{https://stacks.math.columbia.edu}.
%
\bibitem{VN} {\sc J. Van Name}, {\sl 
Ultraparacompactness and Ultranormality}, preprint (2013), \url{https://arxiv.org/abs/1306.6086}.
%
\bibitem{W} {\sc A.  Weil},  "Adeles and algebraic groups",  Progress in Mathematics vol. 23,  Birkh\"auser, 1982.
\end{thebibliography}
\end{document}